\documentclass[dvipsnames,DIV=14]{scrartcl}

\usepackage{amsmath,amssymb,amsthm,bm}
\usepackage[T1]{fontenc}

\usepackage{algorithm,algorithmic}
\usepackage[ruled,vlined,linesnumbered,algo2e,algonl]{algorithm2e}
\SetAlCapHSkip{0pt}
\DontPrintSemicolon
\SetNlSty{normalfont}{}{:}
\SetKwProg{Fn}{function}{}{end function}

\makeatletter
\let\orig@algocf@nl@sethref\algocf@nl@sethref
\renewcommand{\algocf@nl@sethref}[1]{
  \orig@algocf@nl@sethref{#1}
  \gdef\@currentcounter{AlgoLine}
}
\makeatother

\usepackage{caption}
\usepackage{graphicx}

\usepackage{booktabs}
\usepackage{multirow}

\usepackage{nicematrix}
\usepackage{authblk}

\usepackage[colorlinks,
linkcolor=MidnightBlue,
citecolor=ForestGreen,
urlcolor=BrickRed]{hyperref}
\usepackage{aliascnt}

\usepackage[noabbrev]{cleveref}
\crefalias{algocf}{algorithm}
\crefalias{AlgoLine}{line}

\graphicspath{{figs/}}

\newenvironment{keywords}{\medskip\textbf{Keywords:}}{}

\newenvironment{MSCcodes}{\medskip\textbf{AMS subject classifications (2020).}}{}

\theoremstyle{plain}
\newtheorem{theorem}{Theorem}

\newaliascnt{lemma}{theorem}
\newtheorem{lemma}[lemma]{Lemma}
\aliascntresetthe{lemma}

\newaliascnt{proposition}{theorem}

\aliascntresetthe{proposition}

\newaliascnt{corollary}{theorem}

\aliascntresetthe{corollary}

\theoremstyle{definition}

\theoremstyle{remark}
\newtheorem{remark}{Remark}

\definecolor{bkgndcolor}{rgb}{1,1,1}
\definecolor{bkgndcolor}{rgb}{0.78,0.93,0.8}

\usepackage{tcolorbox}
\tcbuselibrary{skins}
\tcbuselibrary{breakable}
\tcbset{shield externalize}

\DeclareMathOperator{\cond}{\kappa}

\DeclareMathOperator{\diag}{diag}

\DeclareMathOperator{\trace}{tr}

\newcommand*{\trans}{^{\top}}
\newcommand*{\transinv}{^{-\top}}

\newcommand*{\wh}{\widehat}

\newcommand{\R}{\mathbb{R}}

\def\Lambda{\varLambda}
\def\Sigma{\varSigma}

\newcommand*{\abs}[1]{\lvert#1\rvert}
\newcommand*{\norm}[1]{\lVert#1\rVert}
\newcommand*{\normb}[1]{\big\lVert#1\big\rVert}

\def\adots{\mathinner{\mkern2mu\raise1pt\hbox{.}\mkern2mu
    \raise4pt\hbox{.}\mkern2mu\raise7pt\hbox{.}\mkern1mu}}

\newcommand*{\ul}{\ensuremath{u_\ell}}

\newenvironment{mf}{
  \list{}{\leftmargin0.2cm\rightmargin\leftmargin}\item\relax
  \begin{tcolorbox}[breakable, 
    colback=Dandelion!10!white, sharp corners, boxrule=0pt,frame hidden, 
    enhanced, borderline west={1mm}{-2mm}{Green}]
    \small\sffamily{\color{Green}\textbf{Max:}}}
  {\end{tcolorbox}\endlist}

\newenvironment{yuqi}{
  \list{}{\leftmargin0.2cm\rightmargin\leftmargin}\item\relax
  \begin{tcolorbox}[breakable, 
    colback=Blue!10!white, sharp corners, boxrule=0pt,frame hidden, 
    enhanced, borderline west={1mm}{-2mm}{RoyalBlue}]
    \small\sffamily{\color{RoyalBlue}\textbf{Yuqi's remark:}}}
  {\end{tcolorbox}\endlist}

\def\printcomment{1}
\ifx\printcomment\undefined
  \usepackage{environ}
  \NewEnviron{killcontents}{}

\fi

\usepackage[normalem]{ulem}

\usepackage{amsmath}
\usepackage{pgfplots}
\usepackage{pgfplotstable}
\usepgfplotslibrary{groupplots}
\usetikzlibrary{calc,patterns,patterns.meta}
\pgfplotsset{compat=1.18}

\pgfplotscreateplotcyclelist{cqrIterations}{
  solid,color=Blue,mark=square,mark size=1.5pt\\
  solid,color=Cyan,mark=o,mark size=1.5pt\\
  solid,color=BurntOrange,mark=asterisk,mark size=2.0pt\\
  solid,color=BrickRed,mark=diamond,mark size=1.8pt\\
  solid,color=Purple,mark=x,mark size=1.9pt\\
  solid,color=Sepia,mark=+,mark options={rotate=-90},mark size=1.9pt\\
}

\pgfplotscreateplotcyclelist{cqrHistories}{
  solid,color=Peach,mark=*,mark size=1.9pt\\
  dash dot,color=Orange!70!Black,mark options={solid},mark=asterisk,mark size=1.9pt\\
  loosely dotted,color=Bittersweet!80!White,mark options={solid},mark=+,mark size=2.1pt\\
  loosely dashed,color=Brown,mark=diamond,mark options={rotate=0,solid},mark size=1.8pt\\
  loosely dash dot,color=Cyan!80!White,mark=square*,mark options={solid,rotate=0},mark size=2.0pt\\
  dashed,color=Aquamarine!80!White,mark=square*,mark options={solid,rotate=-45},mark size=1.9pt\\
  densely dash dot,color=Brown,mark=o,mark options={solid,rotate=0},mark size=1.8pt\\
  densely dashed,color=CarnationPink,mark=triangle,mark options={solid,rotate=0},mark size=1.8pt\\
  dash pattern={on 4pt off 1pt on 1pt off 1pt},color=Red,mark=triangle,mark options={solid,rotate=90},mark size=1.8pt\\
  densely dotted,color=Blue,mark=triangle,mark options={solid,rotate=180},mark size=1.8pt\\
  dotted,color=MidnightBlue,mark=triangle,mark options={solid,rotate=270},mark size=1.8pt\\
}
\pgfplotscreateplotcyclelist{cqrTimingBars}{
  fill=Cyan!40,postaction={pattern=north east lines,pattern color=black!99}\\
  fill=SeaGreen!40,postaction={pattern=north west lines,pattern color=black!99}\\
  fill=OliveGreen,postaction={pattern=crosshatch,pattern color=white}\\
  fill=Violet,postaction={pattern=north east lines,pattern color=white}\\
  fill=Sepia,postaction={pattern=north west lines,pattern color=white}\\
  fill=Apricot!60,postaction={pattern=crosshatch,pattern color=black!99}\\
  fill=LimeGreen,postaction={pattern=grid,pattern color=white}\\
  fill=Black,postaction={pattern=dots,pattern color=white}\\
}

\pgfplotsset{
  every axis/.append style={
    tick align=inside,
  },
  every axis plot post/.append style={
    line width=1pt
  },
  every axis legend/.append style={
    line width=.8pt,
    draw=Black!60,
    fill=white,
    legend cell align=left,
    column sep=3pt,
    inner sep=3pt,
    font=\footnotesize
  },
  cqr/common/.style={scaled ticks=false},
  cqr/log-grid/.style={grid=both,
    major grid style={draw=black!20,line width=.28pt},
    minor grid style={draw=black!14,densely dotted,line width=.24pt}},
  cqr/iterations/.style={
    cqr/common,cqr/log-grid,cycle list name=cqrIterations,
    scale only axis=true,
    width=\dimexpr(\linewidth-3.15cm)/5\relax,
    height=1.7cm,
    xmode=log,ymode=log,log basis x=10,log basis y=10,
    xmin=1,xmax=1e15,ymin=1,ymax=1e15,
    xtick={1,1e3,1e6,1e9,1e12,1e15},
    ytick={1,1e3,1e6,1e9,1e12,1e15},
    tick label style={font=\scriptsize},
    xticklabel style={rotate=45, anchor=north, xshift=-5pt},
    title style={font=\small,yshift=-0.08cm},
    minor tick num=2,
    every axis plot/.append style={line width=.65pt,solid},
    clip mode=individual},
  cqr/iteration-legend/.style={
    legend columns=3},
  cqr/histories/.style={
    cqr/common,cycle list name=cqrHistories,
    scale only axis=true,
    width=\dimexpr(\linewidth-4cm)/3\relax,
    height=2.5cm,
    xmin=0,xmax=6,xtick={0,1,2,3,4,5,6},
    ymin=1,ymax=1e15,ymode=log,log basis y=10,
    ytick={1,1e5,1e10,1e15},
    grid=both,major grid style={draw=black!17,line width=.3pt},
    minor grid style={draw=black!12,densely dotted,line width=.23pt},
    minor tick num=1,tick label style={font=\small},
    xlabel={Iterations},ylabel={},
    title style={font=\normalsize},
    every axis plot/.append style={line width=.72pt,solid},
    unbounded coords=jump},
  cqr/history-legend/.style={
    legend columns=3},
  cqr/timings/.style={
    cqr/common,cycle list name=cqrTimingBars,
    scale only axis=true,
    width=\dimexpr(\linewidth-6cm)\relax,
    height=5cm,
    ymin=0,xmin=.5,xmax=10.5,xtick={1,2,...,10},
    xticklabels={SCQR,mpSCQR,mpEVDQR,LUCQR,mpLUCQR,mpHQR,
      RCQR-ss,mpRCQR-ss,RCQR-g,mpRCQR-g},
    x tick label style={font=\footnotesize,rotate=24,anchor=north east},
    y tick label style={font=\normalsize},ylabel={Time (ms)},
    title style={font=\normalsize},
    grid=major,major grid style={draw=black!16,line width=.3pt},
    bar width=16pt,enlarge x limits=false,
    legend style={at={(1.02,.91)},anchor=north west,inner sep=4pt,row sep=1pt},
    legend cell align=left}
}

\def\FirstWord#1 #2\Stop{#1}
\newcommand{\AddTimingSeries}[1]{
  \addplot+[ybar stacked,mark=none,draw=black!40] table[x=index,y=Gram]{#1};
  \addlegendentry{$X^\top X$}
  \addplot+[ybar stacked,mark=none,draw=black!40] table[x=index,y=LGram]{#1};
  \addlegendentry{$L^\top L$}
  \addplot+[ybar stacked,mark=none,draw=black!40] table[x=index,y=Sketch]{#1};
  \addlegendentry{Sketching}
  \addplot+[ybar stacked,mark=none,draw=black!40] table[x=index,y=FP32]{#1};
  \addlegendentry{Binary32 factorizations}
  \addplot+[ybar stacked,mark=none,draw=black!40] table[x=index,y=FP64]{#1};
  \addlegendentry{Binary64 factorizations}
  \addplot+[ybar stacked,mark=none,draw=black!40] table[x=index,y=Solve]{#1};
  \addlegendentry{$XR_{\mathrm{p}}^{-1}$}
  \addplot+[ybar stacked,mark=none,draw=black!40] table[x=index,y=Accumulate]{#1};
  \addlegendentry{Accumulate $R_{\mathrm{p}}$}
  \addplot+[ybar stacked,mark=none,draw=black!40] table[x=index,y=Other]{#1};
  \addlegendentry{Other}
  \addplot[draw=none,only marks,mark=none,forget plot,
  point meta=explicit symbolic,
  nodes near coords={
    \expandafter\FirstWord\pgfplotspointmeta\Stop
  },
  nodes near coords style={
    font=\scriptsize\bfseries,anchor=south,yshift=0pt
  }]
  table[x=index,y=Total,meta=IterLabel]{#1};
}

\newif\ifsimax
\simaxfalse

\title{Three mixed-precision preconditioners\\for CholeskyQR}

\author[1]{Yuqi Liu}
\author[1]{Massimiliano Fasi}
\affil[1]{School of Computer Science,
University of Leeds, Woodhouse Lane,
Leeds LS2 9JT, UK}
\date{\today}

\begin{document}

\maketitle

\begin{abstract}
  CholeskyQR is an orthogonalization algorithm particularly well suited to
  modern computing architectures, but it may fail to produce numerically
  orthonormal columns when applied to ill-conditioned matrices.
  A customary solution is to first reduce the condition number of the input
  matrix by means of a preconditioner, and then apply CholeskyQR to the
  preconditioned matrix. We propose three mixed-precision techniques for
  constructing such preconditioners: mpSCQR, based on a shifted Cholesky
  factorization; mpEVDQR, based on the eigendecomposition; and mpLUCQR, based on
  a pivoted LU factorization. We prove bounds on the condition number of the
  matrix after mpSCQR, mpLUCQR, or a shifted variant of mpEVDQR has been
  applied.
  Our experiments compare the proposed preconditioners with existing
  deterministic and randomized alternatives. The results suggest that mpEVDQR
  and mpLUCQR can be substantially faster than traditional shiftedCholeskyQR,
  with mpLUCQR attaining performance comparable to that of the most efficient
  randomized preconditioners.
\end{abstract}



\section{CholeskyQR and preconditioning}
\label{sec:introduction}

Let $X\in\R^{m\times n}$ be a matrix of full column rank, with $m\geq n$.
A (column) orthogonalization of $X$ is a factorization that satisfies
\begin{equation}
\label{eq:qr}
X=QR, \qquad Q\trans Q=I_n, \qquad Q\in\R^{m\times n}, \qquad R\in\R^{n\times n}.
\end{equation}
When $R$ is upper triangular, \eqref{eq:qr} is a QR factorization.

Many numerical algorithms require orthogonalizing matrices for which
\mbox{$m \gg n$}. This is the case, for example, when constructing orthonormal
bases in Krylov subspace methods such as GMRES~\cite{sasc86} and the Lanczos
method~\cite{lanc50}, when maintaining the search spaces in eigensolvers such as
PINVIT~\cite{neym02} and LOBPCG~\cite{dsyg18}, and when computing approximate
ranges in randomized low-rank approximation~\cite{hmt11}.

Many algorithms are available for computing~\eqref{eq:qr}, and \emph{Householder
  QR} (HQR) is arguably the most widely used~\cite[sect.~9]{hima22}. HQR is
unconditionally stable, in that the relative backward error and the loss of
orthogonality of its computed QR factors are, up to constants that depend on $m$
and $n$, of order the unit roundoff $u$. The orthogonalization is obtained by
applying to $X$ a sequence of Householder transformations, and requires
$2mn^{2} -2n^{3}/3$ floating-point operations (flops), if $Q$ is only needed as
a product of transformations~\cite[p.~249]{gova13}, and $4mn^{2} - 4n^3/3$
flops, if $Q$ is needed explicitly~\cite[p.~255]{gova13}. The dependencies
require frequent synchronizations and communication that can limit parallel
performance.

\begin{algorithm2e}[t]
  \caption{CholeskyQR}
  \label{alg:choleskyqr}
  \KwIn{$X \in \R^{m\times n}$ with $m \ge n$.}
  \KwOut{$Q$ and $R$ as in~\eqref{eq:qr}.}
  \Fn{$[Q,R]=\texttt{CHOLESKYQR}(X)$}{
    Compute $G\gets X\trans X$.\;\label{ln:cholqr_step1}
    Compute the Cholesky factorization $R\trans R=G$.\;\label{ln:cholqr_step2}
    Compute $Q\gets XR^{-1}$.\;\label{ln:cholqr_step3}
  }
\end{algorithm2e}

On modern architectures, better performance can typically be achieved by using
algorithms with lower communication cost, such as CholeskyQR~\cite{fnyy14} in
\cref{alg:choleskyqr}. The algorithm yields a QR factorization because $R$ is
upper triangular and
\begin{equation*}
  Q\trans Q=R\transinv X\trans XR^{-1}=I_n,
\end{equation*}
where $I_{n}$ is the identity matrix of order $n$.

Forming the symmetric Gram matrix on \cref{ln:cholqr_step1} requires $mn^{2}$
flops, computing its Cholesky factorization on \cref{ln:cholqr_step2} requires
$n^{3}/3$ flops, and solving the linear system on \cref{ln:cholqr_step3}
requires $mn^{2}$ flops. Therefore, CholeskyQR requires $2mn^{2} + n^{3}/3$
flops, which is comparable with the flop counts of the fastest orthogonalization
methods available, such as classical and modified
Gram--Schmidt~\cite[pp.~254--255]{gova13}. All three steps in CholeskyQR are
rich in level-3 BLAS operations, and they can be naturally parallelized with
only one global reduction when forming the Gram matrix. These properties make
CholeskyQR particularly attractive on multicore processors, distributed-memory
systems, and GPUs~\cite{ytd15}.

Since $\cond_{2}(G) = \cond_2(X\trans X)=\cond_2(X)^2$, the loss of
orthogonality of CholeskyQR depends quadratically on $\cond_{2}(X)$. A natural
solution is to apply CholeskyQR twice. Yamamoto et al.~\cite{ynyf15} show that
this variant, known as CholeskyQR2, computes a matrix $\wh Q$ such that
\begin{equation*}
  \norm{\wh Q\trans\wh Q-I_n}_2\leq 6\bigl(mn+n(n+1)\bigr)u,
\end{equation*}
provided that $X$ satisfies
\begin{equation}
\label{eq:cqr2_condition}
8\cond_2(X)\sqrt{mnu+n(n+1)u}\leq 1.
\end{equation}

Condition~\eqref{eq:cqr2_condition} is certainly violated if
$\cond_2(X)>u^{-1/2}$, and CholeskyQR2 is no longer guaranteed to produce a
numerically orthogonal matrix in finite precision. Moreover, the computed Gram
matrix may be indefinite, leading to a breakdown in the Cholesky factorization.
In such cases, one can reduce the condition number by using
\emph{shiftedCholeskyQR} (SCQR), which uses \cref{alg:choleskyqr} with
\cref{ln:cholqr_step2} replaced by
\begin{equation*}
R_{\mathrm{p}}\trans R_{\mathrm{p}}=G+sI_n,
\end{equation*}
for a small shift $s>0$.
Fukaya et al.~\cite{fkny20} show that if the shift is chosen as
\begin{equation}
  \label{eq:shift-fknyy}
s=11\bigl(mn+n(n+1)\bigr)u\norm{X}_2^2
\end{equation}
and the matrix $X$ satisfies
\begin{equation}
\label{eq:SCQR_assum}
\cond_2(X)\leq\frac{u^{-1}}{96\bigl(mn+n(n+1)\bigr)},
\end{equation}
then the Cholesky factorization of $G + sI_n$ does not break down and the
computed preconditioned matrix $\wh Q_{\mathrm{p}} = X \wh R_{\mathrm{p}}^{-1}$
satisfies
\begin{equation}
\label{eq:SCQR_ortho}
\cond_2(\wh Q_{\mathrm{p}})\leq\frac{1}{8\sqrt{\bigl(mn+n(n+1)\bigr)u}}.
\end{equation}
Condition~\eqref{eq:SCQR_assum} shows that SCQR can be applied to input matrices
with condition numbers as large as $u^{-1}$, while~\eqref{eq:SCQR_ortho}
guarantees that the preconditioned matrix $\wh Q_{\mathrm{p}}$
satisfies~\eqref{eq:cqr2_condition} and that CholeskyQR2 can be applied to it.

SCQR can be naturally understood as a form of preconditioning.
A nonsingular matrix $R_{\mathrm{p}}\in\R^{n\times n}$ is a preconditioner for
CholeskyQR if the condition number of the preconditioned matrix
\mbox{$Q_{\mathrm{p}}=XR_{\mathrm{p}}^{-1}$} is substantially lower than that of $X$.
Ideally, $\cond_{2}(Q_{\mathrm{p}})$  should be sufficiently small that
CholeskyQR or CholeskyQR2 can be applied reliably to $Q_{\mathrm{p}}$,
yielding $Q_{\mathrm{p}}=QR_{\mathrm{c}}$ and an orthogonalization
$X=Q(R_{\mathrm{c}}R_{\mathrm{p}})$.

Other preconditioning strategies have been proposed in the literature. Terao,
Ozaki, and Ogita~\cite{too20} propose \emph{LU--CholeskyQR} (LUCQR), which first
computes a pivoted LU factorization \mbox{$PX=LU$}, then forms the Gram matrix
\mbox{$G=L\trans L$} and uses its Cholesky factorization $S\trans S=G$ to
compute the preconditioner $R_{\mathrm{p}}=SU$. In exact arithmetic, we have
\begin{equation*}
X\trans X=U\trans L\trans PP\trans LU
=U\trans S\trans SU=R_{\mathrm{p}}\trans R_{\mathrm{p}},
\end{equation*}
which shows that  $XR_{\mathrm{p}}^{-1}$ has orthonormal columns.

It can be shown~\cite{too20} that if the computed LU decomposition is such that
\begin{align}
8\cond_2(\wh L)\sqrt{mnu+n(n+1)u}&\leq 1,\label{eq:lucqr_condition1}\\
64n^2u\cond_2(\wh L)\cond_2(\wh U)&\leq 1,\label{eq:lucqr_condition2}
\end{align}
then one CholeskyQR step will produce
an orthogonality error and a relative residual of order $u$.
Conditions~\eqref{eq:lucqr_condition1} and~\eqref{eq:lucqr_condition2}
only feature the condition numbers of the LU factors of $X$, not the
condition number of $X$ itself. Therefore, unlike SCQR, LUCQR is not
guaranteed to work for all matrices $X$ with condition numbers of
order $u^{-1}$. From the perspective of rounding-error analysis, if
$\cond_2(X)\approx\cond_2(L)$, the sufficient
condition~\eqref{eq:lucqr_condition1} is equivalent to~\eqref{eq:cqr2_condition}.
Thus whenever CholeskyQR2 is at risk of breakdown, so is LU--CholeskyQR.
However, the practical advantage of LUCQR is that $L$ is often much better
conditioned than~$X$ \cite{goto08}, making the Cholesky factorization of
$L\trans L$ safer than that of $X\trans X$.

Our goal is to develop efficient preconditioners for CholeskyQR using mixed
precision. Mixed-precision preconditioning, first proposed by Carson and
Higham~\cite{cahi18}, is now a well established technique for linear systems,
but relatively little attention has been devoted to mixed-precision
preconditioning for CholeskyQR. In the context of eigenvalue problems, Kressner,
Ma, and Shao~\cite{kms23} suggested preconditioning CholeskyQR with the triangular
factor of a low-precision Householder QR factorization.
Their algorithm is summarized in \cref{alg:MPHQR}. Pranesh and Tisseur
considered using as preconditioner a low-precision LU factorization. Their work
was presented at conferences~\cite{tiss21,tipr21}, but, to the best of our knowledge,
a full account is not yet publicly available. Approaches based on randomized sketching
have also been considered in the literature~\cite{bala22,gaip24}, and, more
recently, randomized preconditioners based on low-precision sketching have been
proposed~\cite{gci26,gzwc26}.

We propose three new mixed-precision preconditioners for CholeskyQR: one is
based on a shifted Cholesky factorization, one on an eigendecomposition, and one
on a pivoted LU factorization. Our error analysis quantifies how much a single
application of these three preconditioners can reduce the condition number of
$X$. This is the subject of \cref{sec:algorithms}; the proofs are in
\cref{sec:conv_mpscqr_proof,sec:conv_mplucqr_proof} in the appendix.

The analysis leads us to view CholeskyQR preconditioning as an iterative
refinement process, where preconditioning can be applied repeatedly as needed
until CholeskyQR or CholeskyQR2 can be applied reliably. In
\cref{sec:experiments}, we compare the performance of the three mixed-precision
preconditioners with that of our implementations of several existing approaches.
The numerical results suggest that, in practice, the proposed methods are
comparable with the state of the art.

\section{Three mixed-precision preconditioners for CholeskyQR}
\label{sec:algorithms}

The preconditioning strategies mentioned in \cref{sec:introduction} are as
expensive as one CholeskyQR step. As the whole orthogonalization usually
consists of one preconditioning step followed by one or two applications of
CholeskyQR, computing and applying the preconditioner can account for a large
fraction of the total computational cost. In this section, we investigate how
low precision can be used to make this step more efficient.

We note that for an ill-conditioned $X$, the computed preconditioner
$\wh R_{\mathrm{p}}$ and the preconditioned matrix
$\wh Q_{\mathrm{p}}$, stored in working precision, should satisfy
\begin{align}
\cond_2(\wh Q_{\mathrm{p}})&<\cond_2(X),\label{eq:conddecrease}\\
\frac{\norm{X-\wh Q_{\mathrm{p}}\wh R_{\mathrm{p}}}_2}
{\norm{\wh Q_{\mathrm{p}}\wh R_{\mathrm{p}}}_2}&=O(u),\label{eq:XQRerror}
\end{align}
where~\eqref{eq:conddecrease} expresses the purpose of preconditioning,
and~\eqref{eq:XQRerror} captures the requirement of a working-precision backward
error in the representation of $X$.
This is necessary because, in most applications, orthogonalization is performed
to obtain a basis for the subspace spanned by $X$, and if the orthogonalization
does not preserve the subspace accurately, the overall algorithm may fail to
converge~\cite{slvo96}.

\subsection{An iterative framework for CholeskyQR preconditioning}

The preconditioners we propose in this section can
satisfy~\eqref{eq:conddecrease}, but the use of low precision may limit
the effect of a single preconditioning step, leaving the preconditioned
matrix $\wh Q_{\mathrm{p}}$ too ill-conditioned for CholeskyQR or CholeskyQR2.
This issue can be alleviated by applying a preconditioner \texttt{PREC}
repeatedly, progressively reducing the condition number until it meets the
requirements of CholeskyQR or CholeskyQR2.

This strategy is shown in \cref{alg:iterprec}.
In this section, we present three mixed-precision
preconditioners for CholeskyQR that can serve as \texttt{PREC}.
For each strategy, we analyze how much a single application can
reduce the condition number.

\begin{algorithm2e}[t]
\caption{Iterative framework for CholeskyQR preconditioning}
\label{alg:iterprec}
\Indm
\KwIn{$X\in\mathbb{R}^{m\times n}$;
CholeskyQR preconditioner $\texttt{PREC}$;
condition number tolerance $\tau$;
maximum number of iterations $j_{\max}$.}
\KwOut{CholeskyQR preconditioner $R_{\mathrm{p}}$;
  preconditioned matrix $Q_{\mathrm{p}}=XR_{\mathrm{p}}^{-1}$.}
\Indp
Initialize $Q_{\mathrm{p}}\gets X$ and $R_{\mathrm{p}}\gets I_n$.\;
\For{$j=1,\dotsc,j_{\max}$}{
    Perform preconditioning
    $[Q_{\mathrm{p}},\widetilde R_{\mathrm{p}}]=\texttt{PREC}(Q_{\mathrm{p}})$.\;
    Compute $R_{\mathrm{p}}\gets \widetilde R_{\mathrm{p}}R_{\mathrm{p}}$.\;
    \If{$\cond_2(Q_{\mathrm{p}})<\tau$}{
        \textbf{break}\;
    }
}
\end{algorithm2e}

\subsection{Mixed-precision shiftedCholeskyQR (mpSCQR)}

For $m \gg n$, the cost of SCQR is dominated by the formation of $G=X\trans X$
and by the solution of the linear system
$Q_{\mathrm{p}} = XR_{\mathrm{p}}^{-1}$. $Q_{\mathrm{p}}$ must be computed in
working precision to satisfy \eqref{eq:XQRerror}, but we can form $G$ in low
precision, obtaining mpSCQR in \cref{alg:MPSCQR}.

\begin{algorithm2e}[t]
\caption{Mixed-precision shiftedCholeskyQR (mpSCQR)}
\label{alg:MPSCQR}
\Indm \KwIn{$X\in\mathbb{R}^{m\times n}$; shift $s>0$; unit roundoffs
  $u\le \ul$.} \KwOut{CholeskyQR preconditioner $R_{\mathrm{p}}$; preconditioned
  matrix $Q_{\mathrm{p}}=XR_{\mathrm{p}}^{-1}$.} \Indp
\Fn{$[Q_{\mathrm{p}},R_{\mathrm{p}}]=\texttt{MPSCQR}(X)$}{ Compute
  $G\gets X\trans X$ in precision $\ul$.\; Compute the Cholesky factorization
  $R_{\mathrm{p}}\trans R_{\mathrm{p}}=G+sI_n$ in precision $u$.\; Compute
  $Q_{\mathrm{p}}\gets XR_{\mathrm{p}}^{-1}$ in precision
  $u$.\;\label{alg-step:mpscqr_XR} }
\end{algorithm2e}

The following theorem bounds the condition number after one mpSCQR step.

\begin{theorem}[Condition number reduction of mpSCQR]
\label{thm:mpscqr-convergence}
Let $X \in \R^{m \times n}$ satisfy
\begin{align}
6n^2u\cond_2(X)&<1,\label{eq:mpscqr-assump-cond}\\
11\bigl((m+3)n\ul+n(n+3)u\bigr)\norm{X}_2^2\le s
&\le \frac{1}{100}\norm{X}_2^2.\label{eq:mpscqr-shift-range}
\end{align}
Then mpSCQR does not break down, and the preconditioned $\wh Q_{\mathrm{p}}$
satisfies
\begin{equation}
\label{eq:mpscqr-main-bound}
\cond_2(\wh Q_{\mathrm{p}})
\le1.4\sqrt{\frac{1+1.1s\cond_2(X)^2/\norm{X}_2^2}{1+0.9s/\norm{X}_2^2}}.
\end{equation}
In particular, if $\sqrt{(m+3)n\ul+n(n+3)u}\cond_2(X)\gg1$ and $s$ is chosen to be
the lower bound in~\eqref{eq:mpscqr-shift-range},
then~\eqref{eq:mpscqr-main-bound} implies that one mpSCQR step reduces the
condition number of $X$ by a factor of
\begin{equation}
\label{eq:mpscqr-factor}
O\bigl(\sqrt{(m+3)n\ul+n(n+3)u}\bigr).
\end{equation}
\end{theorem}

\begin{proof}
See \cref{sec:conv_mpscqr_proof}.
\end{proof}

Assumption~\eqref{eq:mpscqr-shift-range} restricts the low-precision unit
roundoff $\ul$, as it implies that
\begin{equation}
\label{eq:mpscqr-small-dim}
(m+3)n\ul+n(n+3)u\le\frac{1}{1100}.
\end{equation}

When $m\gg n$ and $\ul\gg u$, the factor in~\eqref{eq:mpscqr-factor} can be
simplified to $O\bigl(\sqrt{mn\ul}\bigr)$, which becomes
$O\bigl(\sqrt{\ul}\bigr)$ if the dependence on the matrix dimensions is
disregarded.

Choosing an appropriate shift can be difficult. The
bound~\eqref{eq:mpscqr-main-bound} shows that convergence is faster the smaller
$s$ is, and the lower bound in~\eqref{eq:mpscqr-shift-range} is too large to
guarantee reasonably fast convergence. Therefore, in practice, we often choose a
smaller $s$, accepting the risk that Cholesky factorization may break down. In
\cref{sub-sec:exp_size}, we show that for mpSCQR this is even more likely to
occur than for SCQR.

\subsection{Mixed-precision eigendecomposition QR (mpEVDQR)}

mpSCQR can be interpreted as correcting the spectrum of the low-precision Gram
matrix~$\wh G$ by adding a uniform shift. A more general approach, suggested by
this observation, is to compute the spectral decomposition of $\wh G$ explicitly
and correct its eigenvalues. This approach is summarized in \cref{alg:MPEVDQR}.
\begin{algorithm2e}[t]
\LinesNumbered
\DontPrintSemicolon
\caption{Mixed-precision eigendecomposition QR (mpEVDQR)}
\label{alg:MPEVDQR}
\Indm
\KwIn{$X\in\mathbb{R}^{m\times n}$;
  $f(\cdot):\R\rightarrow\R$;
  unit roundoffs $u\le \ul$.}
\KwOut{CholeskyQR preconditioner $R_{\mathrm{p}}$; preconditioned matrix
  $Q_{\mathrm{p}}=XR_{\mathrm{p}}^{-1}$.}
\Indp
\Fn{$[Q_{\mathrm{p}},R_{\mathrm{p}}]=\texttt{MPEVDQR}(X)$}{
    Compute $G\gets X\trans X$ in precision $\ul$.\;
    Compute an eigendecomposition
    $G=V\Lambda V\trans$ in precision $u$.\;
    Correct the spectrum and compute
    $S_\Lambda\gets f(\Lambda)^{1/2}$ in precision $u$.\;
    \label{alg-step:mpevdqr_speccorrect}
    Compute $R_{\mathrm{p}}\gets S_\Lambda V\trans$
    in precision $u$.\;
    \label{alg-step:mpevdqr_genR}
    Compute $Q_{\mathrm{p}}\gets XV S_\Lambda^{-1}$
    in precision $u$.\;
    \label{alg-step:mpevdqr_XR}
}
\end{algorithm2e}

Once the eigenvalues are available, they can be corrected using any function
$f(\cdot):\R\rightarrow\R$ satisfying $f(\wh\lambda)>0$ for every computed
eigenvalue $\wh\lambda$ of $\wh G$. The uniform shift used in SCQR and mpSCQR
can be understood as a special case with $f(x)=x+s$. The main drawback of this
approach is that the preconditioner $R_{\mathrm{p}}=S_\Lambda V\trans$ returned
by mpEVDQR is generally dense rather than triangular, and this has two negative
consequences: first, updating $X$ on \cref{alg-step:mpevdqr_XR} requires a
matrix--matrix product, which is more expensive than the triangular solve
required in mpSCQR; second, mpEVDQR followed by CholeskyQR does not produce a
standard QR factorization. The latter issue can be addressed by performing, on
\cref{alg-step:mpevdqr_genR}, an additional QR factorization of
$S_\Lambda V\trans$, whose cost will be negligible if $m \gg n$.

The convergence behavior of \cref{alg:MPEVDQR} depends on the choice of the
correction function $f(\cdot)$. If the eigendecomposition satisfies
\begin{equation*}
\norm{\wh V\wh\Lambda\wh V\trans-\wh G}_2
\leq c_1n^2u\norm{\wh G}_2,
\qquad
\norm{\wh V\wh V\trans-I_n}_2
\leq c_2n^2u,
\end{equation*}
for some moderate constants $c_1$ and $c_2$, then setting $f(x)=x+s$ with the
smallest $s$ allowed by~\eqref{eq:mpscqr-shift-range} will achieve a reduction
rate similar to that of mpSCQR, which is $O(\sqrt{\ul})$ if the dependence on
the matrix dimensions is disregarded. In practice, we use the absolute value
correction
\begin{equation}
  \label{eq:eig-corr}
f(x)=\abs{x}.
\end{equation}
This correction produces positive values unless one of the computed eigenvalues
is exactly zero. In \cref{sub-sec:exp_size}, we show experimentally that the
reduction rate of mpEVDQR with the absolute value correction~\eqref{eq:eig-corr}
is $O(\sqrt{\ul})$.

\subsection{Mixed-precision LU--CholeskyQR (mpLUCQR)}
The lower triangular factor $L$ of the partially pivoted LU factorization
$PX=LU$ is often observed to be much better conditioned than $X$. Therefore, we
may hope that Cholesky factorization of $G = L\trans L$ will succeed, even if
the LU decomposition and the Gram matrix are computed in low precision. In fact,
if $\cond_2(L)\ll\cond_2(X)$, then $\wh G$ is substantially more likely to
remain positive definite and admit a successful Cholesky factorization, since
\begin{equation*}
\cond_2(L\trans L)=\cond_2(L)^2
\ll\cond_2(X)^2=\cond_2(X\trans X).
\end{equation*}
The resulting Cholesky factor is then combined with $U$ to construct the
preconditioner. This is the mpLUCQR algorithm, summarized in \cref{alg:MPLUCQR}.
\begin{algorithm2e}[t]
\caption{Mixed-precision LU--CholeskyQR (mpLUCQR)}
\label{alg:MPLUCQR}
\Indm
\KwIn{$X\in\mathbb{R}^{m\times n}$; unit roundoffs $u\le \ul$.}
\KwOut{CholeskyQR preconditioner $R_{\mathrm{p}}$;
preconditioned matrix $Q_{\mathrm{p}}=XR_{\mathrm{p}}^{-1}$.}
\Indp
\Fn{$[Q_{\mathrm{p}},R_{\mathrm{p}}]=\texttt{MPLUCQR}(X)$}{
    Compute the pivoted LU factorization $LU=PX$ in precision $\ul$.\;
    Compute $G\gets L\trans L$ in precision $\ul$.\;
    Compute the Cholesky factorization $S\trans S=G$ in precision $u$.\;
    Compute $R_{\mathrm{p}}\gets SU$ in precision $u$.\;
    Compute $Q_{\mathrm{p}}\gets XR_{\mathrm{p}}^{-1}$ in precision $u$.\;
}
\end{algorithm2e}

The following theorem bounds the condition number after one mpLUCQR step.
\begin{theorem}[Condition number reduction of mpLUCQR]
\label{thm:mplucqr-convergence}
Let the computed LU factors satisfy
\begin{align}
n(n+1)u\cond_2^2(\wh L)
&\le \frac{1}{64},\label{eq:mplucqr-assump-u-L}\\
n^2u\cond_2(\wh L)\cond_2(\wh U)
&\le \frac{1}{64},\label{eq:mplucqr-assump-u-LU}\\
mn\ul\cond_2^2(\wh L)
&\le \frac{1}{64},\label{eq:mplucqr-assump-ul}
\end{align}
and assume that
\begin{equation}
\label{eq:mplucqr-assump-chi}
\frac{\chi_{\max}}{\chi_{\min}}
\le
\min\left\{
1+(n\rho+\sqrt{n})\ul\cond_2(X),
\frac{1}{64n^2u\cond_2(\wh L)\cond_2(\wh U)}
\right\},
\end{equation}
where $\rho = \normb{\abs{\wh L}\abs{\wh U}}_2 / \norm{X}_2$ is the growth
factor of $X$ and
\begin{equation}
  \label{eq:chimin-chimax-def}
\chi_{\min}=\min_{x\neq0}
\frac{\norm{Xx}_2}{\norm{\wh L\wh Ux}_2},\qquad
\chi_{\max}=\max_{x\neq0}
\frac{\norm{Xx}_2}{\norm{\wh L\wh Ux}_2}.
\end{equation}
Then the Cholesky factorization in mpLUCQR does not break down, the computed
triangular factor $\wh R_{\mathrm{p}}$ is nonsingular, and the preconditioned
$\wh Q_{\mathrm{p}}$ satisfies
\begin{equation}
\label{eq:mplucqr-main-bound}
\cond_2(\wh Q_{\mathrm{p}})
\le1.12\cdot\frac{\chi_{\max}}{\chi_{\min}}
\le1.12\bigl(1+(n\rho+\sqrt{n})\ul\cond_2(X)\bigr).
\end{equation}
Consequently, when $(n\rho+\sqrt{n})\ul\cond_2(X)\gg1$,
one mpLUCQR step reduces the condition number of $X$ by a
factor of $O\bigl((n\rho+\sqrt{n})\ul\bigr)$.
\end{theorem}

\begin{proof}
See \cref{sec:conv_mplucqr_proof}.
\end{proof}

However, the lower triangular factor $L$ is not guaranteed to be well
conditioned, and, in particular, it need not be better conditioned than $X$.
For example, the matrix $X$ with
\begin{equation*}
X_{ij}
=
\begin{cases}
1, & i=j,\\
-1, & i>j,\\
0, & i<j,
\end{cases},
\qquad
(X^{-1})_{ij}
=
\begin{cases}
1, & i=j,\\
2^{\,i-j-1}, & i>j,\\
0, & i<j,
\end{cases}
\end{equation*}
has an LU decomposition $PX = LU$ with $L = X$ and $P = U = I_n$.
Since $\norm{L}_{\infty}=n$ and $\norm{L^{-1}}_{\infty}=2^{n-1}$, we have that
\begin{equation*}
\cond_2(X)=\cond_2(L)\geq \frac{\cond_{\infty}(L)}{n}=2^{n-1}.
\end{equation*}

The assumptions in \cref{thm:mplucqr-convergence} show that mpLUCQR can be
ineffective and even break down if $L$ is ill conditioned. When $m\gg n$ and
$\ul\gg u$, condition~\eqref{eq:mplucqr-assump-ul} can become even more
restrictive than~\eqref{eq:mplucqr-assump-u-LU} if
$\cond_2(\wh L) \approx \cond_2(X)$. Thus, the theorem suggests that the
effectiveness of mpLUCQR relies on the commonly observed behavior
\begin{equation*}
\cond_2(\wh L)\ll\cond_2(X).
\end{equation*}

\begin{remark}\label{rmk:mplucqr-assump-chi-1}
Taking into account the rounding error incurred when converting $X$ to
precision $\ul$, the LU factorization with partial pivoting computed
in precision $\ul$ satisfies the componentwise backward-error bound
\begin{equation}
\label{eq:mplucqr-E1}
PX-\wh L\wh U=E_1,
\qquad
\abs{E_1}\le n \ul\abs{\wh L}\abs{\wh U}+\ul P\abs{X}.
\end{equation}
If $(n\rho+\sqrt{n})\ul\cond_2(X)<1$, then \eqref{eq:mplucqr-E1} implies
\begin{equation*}
\frac{\chi_{\max}}{\chi_{\min}}
\le\frac{1+(n\rho+\sqrt{n})\ul\cond_2(X)}{1-(n\rho+\sqrt{n})\ul\cond_2(X)}
\approx1+O\bigl((n\rho+\sqrt{n})\ul\cond_2(X)\bigr).
\end{equation*}
This is consistent with the first term in~\eqref{eq:mplucqr-assump-chi}. In
mixed precision, $(n\rho+\sqrt{n})\ul\cond_2(X)$ can be much larger than one,
so~\eqref{eq:mplucqr-assump-chi} cannot in general be obtained
from~\eqref{eq:mplucqr-E1}. We therefore impose it as a separate assumption to
quantify the quality of the low-precision pivoted LU factorization as a
preconditioner for $X$. Although a general proof for nonsingular square matrices
is still unavailable, such behavior has been widely observed and also forms a
theoretical basis for many mixed-precision linear solvers~\cite{cahi17}.
\end{remark}

\begin{remark}
  The second term in~\eqref{eq:mplucqr-assump-chi} is a rather restrictive
  assumption,
which requires the unrealistic $\chi_{\max}/\chi_{\min}=1$ if the bound
in~\eqref{eq:mplucqr-assump-u-LU} is attained. This restriction is partly caused
by worst-case dimension-dependent error estimates, which are only of theoretical
interest. In practice, mpLUCQR typically remains effective even when this
condition is not satisfied.
\end{remark}

\section{Numerical experiments}
\label{sec:experiments}
In this section, we demonstrate the performance of
\cref{alg:MPSCQR,alg:MPEVDQR,alg:MPLUCQR} via numerical experiments conducted on
an NVIDIA RTX PRO 6000 Blackwell Server Edition GPU with 96 GiB of memory. Our
implementation uses the NVIDIA CUDA Toolkit~13.0 and the cuBLAS~13.5.1 and
cuSOLVER~12.2.2 libraries. The source code is available online.
\footnote{\url{https://github.com/IzayoiYuuki/ThreeMixed-precisionPreconditionersforCholeskyQR_public}.}

In all experiments, we use as working and low precision binary64 and binary32
arithmetic, so that $u = 2^{-53}$ and $\ul = 2^{-24}$. Unless otherwise stated,
we set $s = 10^{-6}$ in \cref{alg:MPSCQR} and use the absolute value
correction~\eqref{eq:eig-corr} in \cref{alg:MPEVDQR}.

The test matrices are generated, for a prescribed condition number $\kappa$, as
\begin{equation}
\label{eq:test-matrix}
\begin{aligned}
X&=U \diag(\sigma_1, \ldots, \sigma_n)V\trans\in\R^{m\times n},\quad
m = 10^{6},\quad \sigma_i=\kappa^{-\frac{i-1}{n-1}},\quad 1\leq i\leq n,
\end{aligned}
\end{equation}
where $U\in\R^{m\times n}$ and $V\in\R^{n\times n}$ are obtained by
orthogonalizing Gaussian matrices.
This matrix is equivalent to mode 3 of the \verb+randsvd+ matrix in the
MATLAB \verb+gallery+.

We compare the new preconditioners against seven existing alternatives:
mono-precision shiftedCholeskyQR (SCQR), mono-precision LU--CholeskyQR (LUCQR),
mixed-precision Householder QR (mpHQR), and mono- and mixed-precision randomized
CholeskyQR (RCQR and mpRCQR, respectively). SCQR is obtained from
\cref{alg:MPSCQR} by setting $u=\ul$ and $s=10^{-14}$, while mpHQR, RCQR, and
mpRCQR are described in \cref{sec:otheralg}. For RCQR and mpRCQR, we test
Gaussian sketching and sparse sign sketching matrices of size $2n \times m$. We
distinguish these two variants using the suffixes ``-g'' and ``-ss'',
respectively. The sparse sign sketch has eight nonzero entries per
column~\cite[p.~456]{matr20}.

In the experiments, we consider two values of $\tau$ for \cref{alg:iterprec}.
The first is $\tau=10$, which reduces the condition number of the preconditioned
matrix to about one. The second is $\tau=10^8$, which represents the cost of
reducing the condition number to a range in which CholeskyQR2 can be applied
successfully. Condition~\eqref{eq:cqr2_condition} depends on $m$ and $n$. Thus,
the choice $\tau = 10^8$ does not guarantee that CholeskyQR2 will be applicable.
In our experience, however, matrices satisfying this threshold can typically be
orthogonalized successfully by CholeskyQR2.

\subsection{Effect of matrix condition number}

We first compare the ten preconditioners on the test matrices
in~\eqref{eq:test-matrix} for $n = 100$ and $\kappa$ between $10^0$ and
$10^{15}$.

\begin{figure}[!tb]
\centering
\begingroup
\newcommand{\IterationSeries}[1]{
  \addplot+ table[x=kappa,y=iter0]{#1};
  \addplot+ table[x=kappa,y=iter1]{#1};
  \addplot+ table[x=kappa,y=iter2]{#1};
  \addplot+ table[x=kappa,y=iter3]{#1};
  \addplot+ table[x=kappa,y=iter4]{#1};
  \addplot+ table[x=kappa,y=iter5]{#1};
}
\begin{tikzpicture}
  \begin{groupplot}[group style={
      group size=5 by 2,
      xticklabels at=edge bottom,
      yticklabels at=edge left,
      horizontal sep=0.35cm,
      vertical sep=1cm
    },
    cqr/iterations]
    \nextgroupplot[title={SCQR},legend to name=FigOneLegend,cqr/iteration-legend]
    \pgfplotstableread[col sep=comma]{figs/fig1_SCQR.csv}\figdata
    \IterationSeries{\figdata}
    \addlegendentry{0 iterations}\addlegendentry{1 iteration}\addlegendentry{2 iterations}
    \addlegendentry{3 iterations}\addlegendentry{4 iterations}\addlegendentry{5 iterations}
    \nextgroupplot[title={mpSCQR}]
    \pgfplotstableread[col sep=comma]{figs/fig1_mpSCQR.csv}\figdata
    \IterationSeries{\figdata}
    \nextgroupplot[title={mpEVDQR}]
    \pgfplotstableread[col sep=comma]{figs/fig1_mpEVDQR.csv}\figdata
    \IterationSeries{\figdata}
    \nextgroupplot[title={LUCQR}]
    \pgfplotstableread[col sep=comma]{figs/fig1_LUCQR.csv}\figdata
    \IterationSeries{\figdata}
    \nextgroupplot[title={mpLUCQR}]
    \pgfplotstableread[col sep=comma]{figs/fig1_mpLUCQR.csv}\figdata
    \IterationSeries{\figdata}
    \nextgroupplot[title={mpHQR}]
    \pgfplotstableread[col sep=comma]{figs/fig1_mpHQR.csv}\figdata
    \IterationSeries{\figdata}
    \nextgroupplot[title={RCQR-ss}]
    \pgfplotstableread[col sep=comma]{figs/fig1_RCQR_Sign.csv}\figdata
    \IterationSeries{\figdata}
    \nextgroupplot[title={mpRCQR-ss}]
    \pgfplotstableread[col sep=comma]{figs/fig1_mpRCQR_Sign.csv}\figdata
    \IterationSeries{\figdata}
    \nextgroupplot[title={RCQR-g}]
    \pgfplotstableread[col sep=comma]{figs/fig1_RCQR_Gaussian.csv}\figdata
    \IterationSeries{\figdata}
    \nextgroupplot[title={mpRCQR-g}]
    \pgfplotstableread[col sep=comma]{figs/fig1_mpRCQR_Gaussian.csv}\figdata
    \IterationSeries{\figdata}
\end{groupplot}
\node[anchor=north] at ([yshift=-0.66cm]group c3r2.south) {\small $\kappa_2(X)$};
\node[anchor=north] at ([yshift=-1.5cm]group c3r2.south){\pgfplotslegendfromname{FigOneLegend}};
\node[rotate=90,anchor=center]
  at ($(group c1r1.north west)!0.5!(group c1r2.south west)
       +(-1cm,0)$)
  {\small $\kappa_2(Q_{\mathrm{p}})$};
\end{tikzpicture}
\endgroup
\vspace{-1.5em}
\caption{Condition number behavior for the ten preconditioning algorithms.
The $x$-axis represents the condition number of $X$, and the $y$-axis represents
the condition number of the preconditioned matrix $Q_{\mathrm{p}}$ after 0 to 5
iterations of~\cref{alg:iterprec}. To ensure that \cref{alg:iterprec} performs
at least five iterations, we set $\tau = 0$. The curve corresponding to zero
iterations represents $Q_{\mathrm{p}}=X$.}
\label{fig:9prec_condhist}
\end{figure}
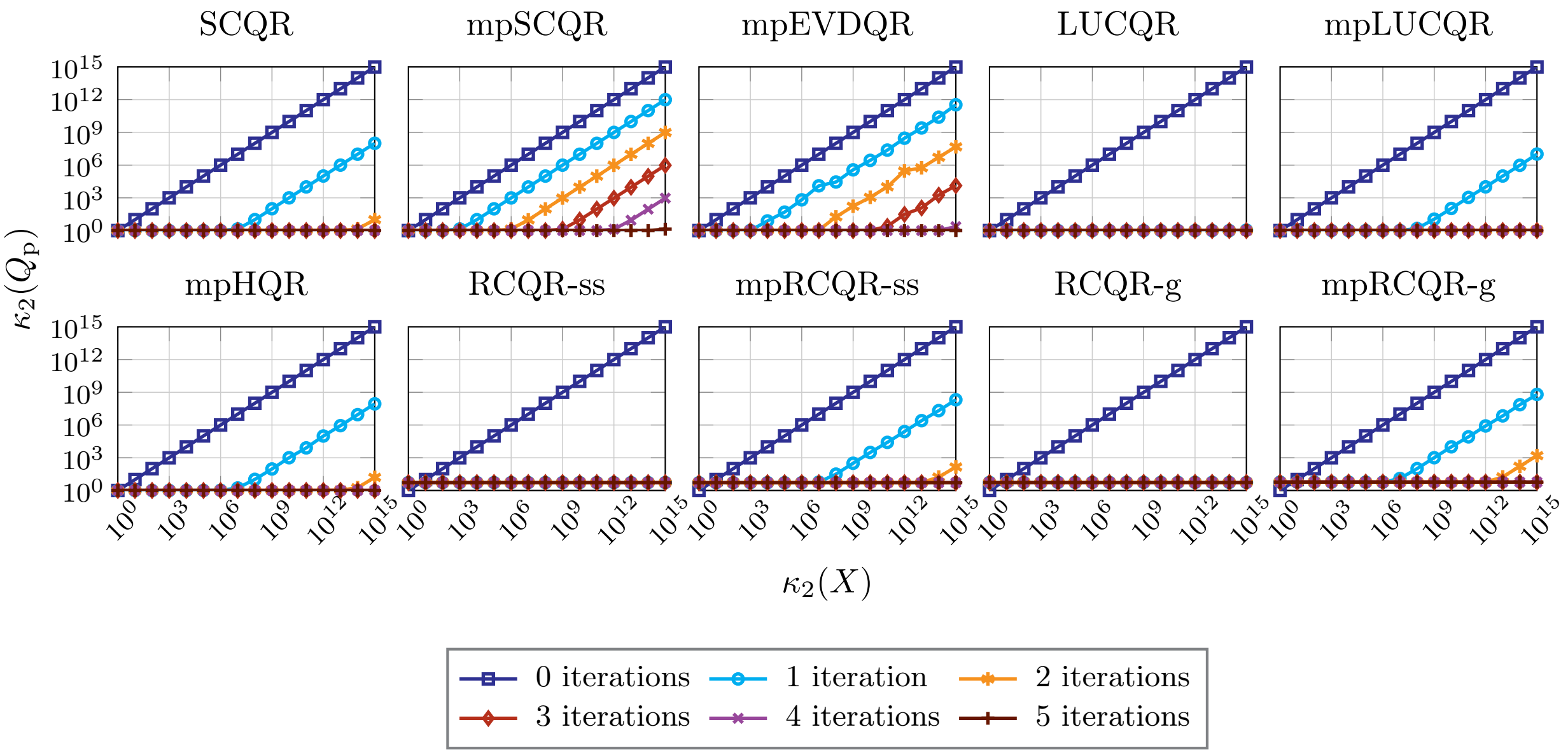

\Cref{fig:9prec_condhist} reports the value of $\cond_2(Q_{\mathrm{p}})$ after one
to five iterations of~\cref{alg:iterprec} for each of the preconditioners.
All ten methods produce a $Q_{\mathrm{p}}$ with condition number of order one within
five iterations.
For each algorithm, the curves corresponding to successive iterations are
approximately parallel, and the vertical gaps between adjacent curves remain
constant as $\cond_2(X)$ increases.

When $\cond_2(X)=10^{15}$, LUCQR, RCQR-g, and RCQR-ss are the most effective
preconditioners, bringing $\cond_2(Q_{\mathrm{p}})$ below $10$ after a single step.
The proposed mpSCQR, mpEVDQR, and mpLUCQR require approximately five, four, and
two iterations, respectively, to bring the condition number down to the same
level. On average, one iteration of mpSCQR and mpEVDQR reduces the condition
number by a factor $10^{-3}$, whereas one iteration of mpLUCQR reduces it by a
factor $10^{-7}$. This agrees with the theoretical factors
discussed in \cref{sec:algorithms}.

\begin{table}[!tb]
\centering
\caption{Classification of the ten preconditioners according to
  whether they (1) use mixed precision and (2) explicitly form $X\trans X$.
  The expected reduction in condition number of a single step is also reported.
  Gaussian and sparse sign sketching variants are grouped for RCQR and mpRCQR.
  }
\begin{tabular}{l@{\qquad}ll@{\quad}l@{\qquad}l@{\quad}l}
\toprule
&& \multicolumn{2}{c}{Form $X\trans X$}
& \multicolumn{2}{c}{Avoid forming $X\trans X$} \\
\midrule
Mono-precision &
  & SCQR & $O(\sqrt{u})$
  & \begin{tabular}[t]{@{}l@{}}
      RCQR \\
      LUCQR
    \end{tabular} & $O(u)$ \\
\midrule
Mixed-precision &
  & \begin{tabular}[t]{@{}l@{}}
      mpSCQR \\
      mpEVDQR
    \end{tabular}
  & $O(\sqrt{\ul})$
  & \begin{tabular}[t]{@{}l@{}}
      mpLUCQR \\
      mpHQR \\
      mpRCQR
    \end{tabular}
  & $O(\ul)$ \\
\bottomrule
\end{tabular}
\label{tab:classification}
\end{table}
The ten preconditioners can be classified according to whether they form
$X\trans X$ explicitly and whether they use more than one precision.
\Cref{tab:classification} lists the algorithms in each class and reports the
expected reduction in condition number of a single step. Taken together,
\cref{tab:classification} and \cref{fig:9prec_condhist} suggest that:
\begin{enumerate}
\item explicitly forming $X\trans X$ introduces a square-root loss in the
  condition number reduction factor, and
\item using lower precision replaces the working-precision unit roundoff $u$ by
  the lower-precision unit roundoff $\ul$.
\end{enumerate}
As we use binary32 as low precision in our experiments, we have that
$\ul\approx\sqrt{u}$ and the use of low precision also corresponds approximately
to a square-root loss.

\subsection{Effect of matrix size}
\label{sub-sec:exp_size}

The factors $O(\sqrt{\ul})$ and $O(\ul)$ derived in \cref{sec:algorithms} were
obtained by disregarding factors that depend on $m$ and $n$. In this section, we
show that the influence of these terms is, to a certain extent, indeed limited.

\begin{figure}[!tb]
\centering
\begingroup
\newcommand{\HistorySeries}[1]{
  \addplot+ table[x=iteration,y=SCQR]{#1};
  \addplot+ table[x=iteration,y=mpSCQR_1e6]{#1};
  \addplot+ table[x=iteration,y=mpSCQR_1e5]{#1};
  \addplot+ table[x=iteration,y=mpEVDQR]{#1};
  \addplot+ table[x=iteration,y=LUCQR]{#1};
  \addplot+ table[x=iteration,y=mpLUCQR]{#1};
  \addplot+ table[x=iteration,y=mpHQR]{#1};
  \addplot+ table[x=iteration,y=RCQR_Sign]{#1};
  \addplot+ table[x=iteration,y=mpRCQR_Sign]{#1};
  \addplot+ table[x=iteration,y=RCQR_Gaussian]{#1};
  \addplot+ table[x=iteration,y=mpRCQR_Gaussian]{#1};
}
\newcommand{\Thresholds}{
  \addplot[draw=black!55,densely dashed,line width=.75pt,forget plot] coordinates {(0,10) (6,10)};
  \addplot[draw=black!55,densely dashed,line width=.75pt,forget plot] coordinates {(0,1e8) (6,1e8)};
}
\begin{tikzpicture}
\begin{groupplot}[group style={group size=3 by 1,horizontal sep=1cm,vertical sep=1cm},cqr/histories]
  \nextgroupplot[title={$n=100$},ylabel={$\kappa_2(Q_{\mathrm{p}})$},legend to name=FigTwoLegend,cqr/history-legend]
    \Thresholds
    \pgfplotstableread[col sep=comma]{figs/fig2_n100.csv}\figdata
    \HistorySeries{\figdata}
    \addlegendentry{SCQR}
    \addlegendentry{mpSCQR ($s=10^{-6}$)}
    \addlegendentry{mpSCQR ($s=10^{-5}$)}
    \addlegendentry{mpEVDQR}
    \addlegendentry{LUCQR}
    \addlegendentry{mpLUCQR}
    \addlegendentry{mpHQR}
    \addlegendentry{RCQR-ss}
    \addlegendentry{mpRCQR-ss}
    \addlegendentry{RCQR-g}
    \addlegendentry{mpRCQR-g}
  \nextgroupplot[title={$n=500$}]
    \Thresholds
    \pgfplotstableread[col sep=comma]{figs/fig2_n500.csv}\figdata
    \HistorySeries{\figdata}
  \nextgroupplot[title={$n=1000$}]
    \Thresholds
    \pgfplotstableread[col sep=comma]{figs/fig2_n1000.csv}\figdata
    \HistorySeries{\figdata}
\end{groupplot}
\node[anchor=north] at ([yshift=-1cm]group c2r1.south)
{\pgfplotslegendfromname{FigTwoLegend}};
\end{tikzpicture}
\endgroup
\caption{Condition number histories of the ten preconditioning algorithms for
  matrices with $\cond_2(X)=10^{15}$ and $n=100$, $500$, and $1000$. For mpSCQR,
  the results with shifts $s=10^{-6}$ and $s=10^{-5}$ are both shown. When
  $s=10^{-6}$ and $n=1000$, the curve for mpSCQR terminates after the first
  iteration because the Cholesky factorization breaks down in the second
  preconditioning step. The dashed horizontal lines indicate the two convergence
  criteria $\cond_2(Q_{\mathrm{p}})=10$ and $\cond_2(Q_{\mathrm{p}})=10^8$.}
\label{fig:curve_condhist}
\end{figure}
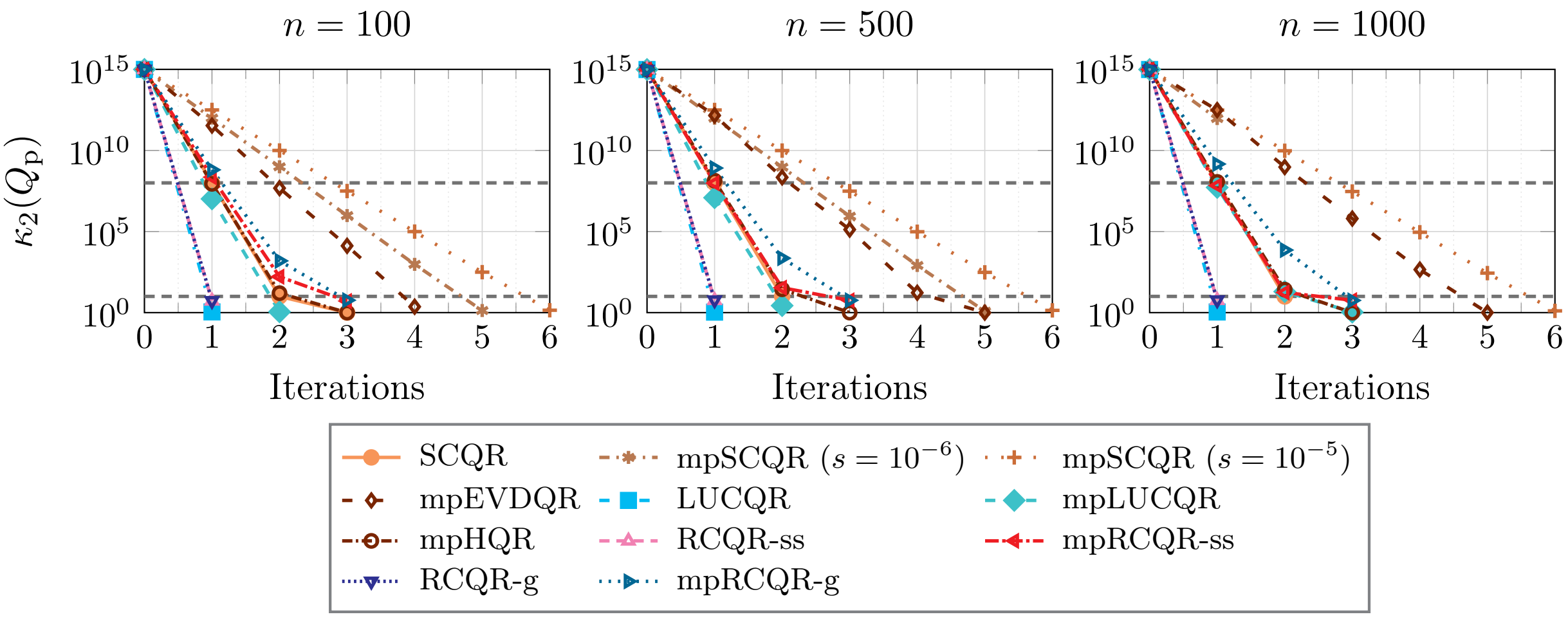

In \cref{fig:curve_condhist}, we test the ten preconditioning algorithms on the
matrices in~\eqref{eq:test-matrix} with $\kappa = 10^{15}$ and $n=100$, $500$,
and $1000$, and we plot the condition number of the preconditioned matrix
$Q_{\mathrm{p}}$ at each iteration. For mpSCQR, we use two different shifts,
$10^{-6}$ and $10^{-5}$, resulting in two separate curves, as indicated in the
legend.

As $n$ increases from $100$ to $1000$, there are no significant changes in the
convergence rates of most algorithms. An exception is mpSCQR with $s=10^{-6}$,
for which the Cholesky factorization breaks down when $n=1000$. We therefore
must use the larger shift $s=10^{-5}$, which slows down convergence. This is a
consequence of the fact that forming the Gram matrix in low precision is
sensitive to the dimension of $X$. mpEVDQR also forms the Gram matrix in low
precision, and its convergence similarly slows down as $n$ increases. However,
the number of iterations until $\cond_2(Q_{\mathrm{p}})<10$ only increases by 25\%,
while $n$ increases by a factor $10$. Therefore, the dependence on $n$ is much
weaker than the bounds in \cref{sec:algorithms} might suggest, and the
approximate reduction factors $O(\sqrt{\ul})$, $O(\sqrt{\ul})$, and $O(\ul)$ are
reasonable descriptions of the observed convergence behavior of mpSCQR, mpEVDQR,
and mpLUCQR, respectively.

\subsection{Performance and profiling}

\begin{figure}[!tb]
\centering
\begingroup
\pgfplotstableread[col sep=comma]{figs/fig3_bars.csv}\cqrbardata
\begin{tikzpicture}
\begin{axis}[
  cqr/timings,
  ymax=200,ytick={0,50,...,200},
  title={Tolerance: $\kappa_2(Q_{\mathrm{p}})<10$}
]
  \AddTimingSeries{\cqrbardata}
\end{axis}
\end{tikzpicture}
\endgroup
\caption{Total runtime required by \cref{alg:iterprec} to achieve
  $\cond_2(Q_{\mathrm{p}})<10$ for a matrix $X$ with $\cond_2(X)=10^{15}$. The
  different sections of each bar represent the individual runtimes of different
  components. The total number of iterations is shown above the bar. The cost of
  estimating $\cond_2(Q_{\mathrm{p}})$ is excluded.}
\label{fig:compbar1e1}
\end{figure}
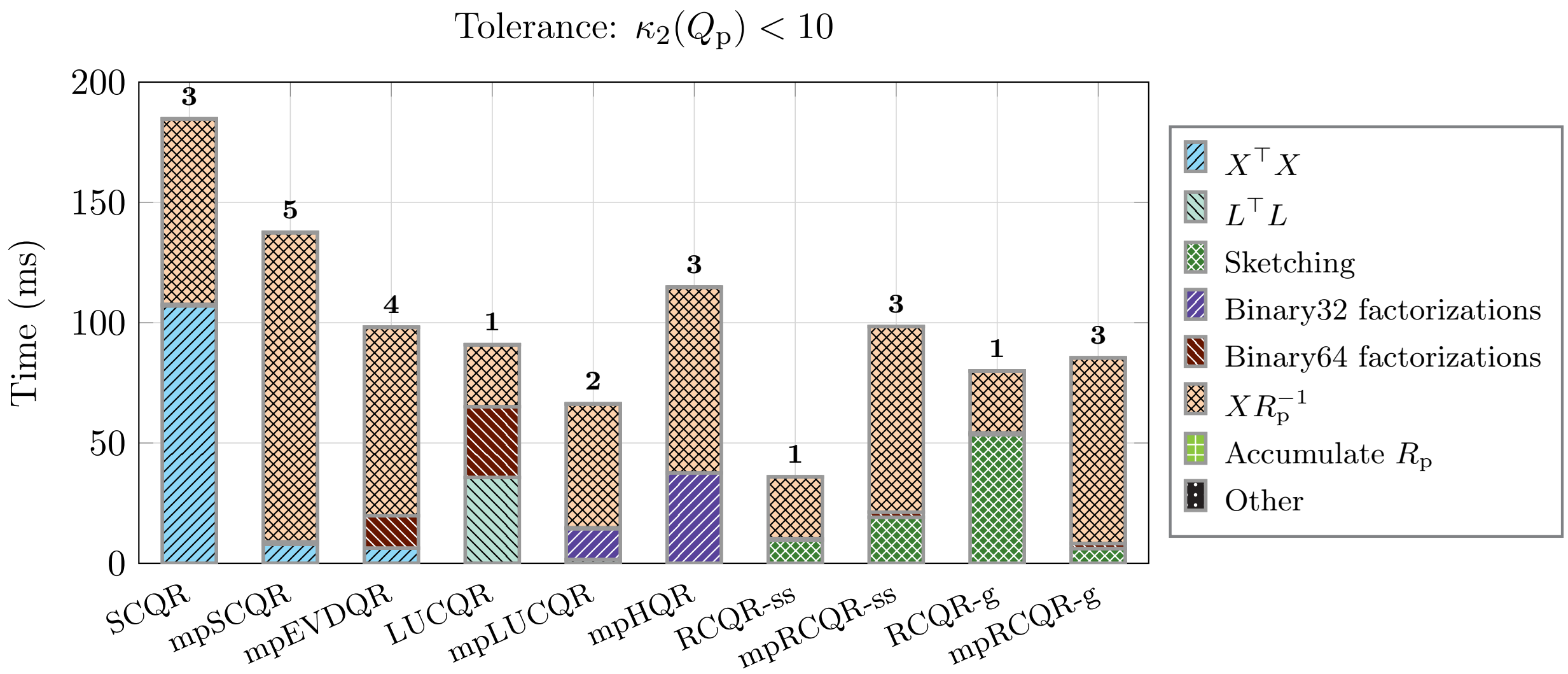

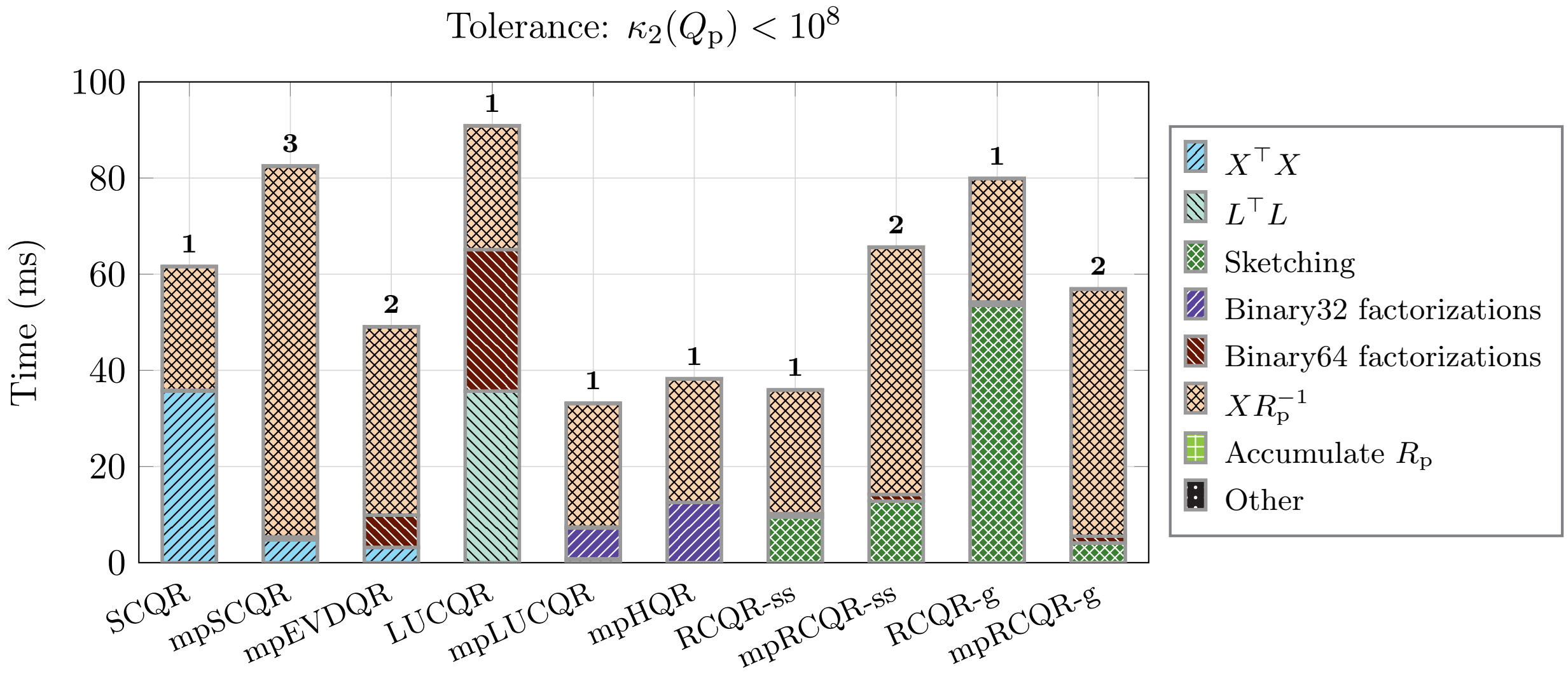
\begin{figure}[!tb]
\centering
\begingroup
\pgfplotstableread[col sep=comma]{figs/fig4_bars.csv}\cqrbardata
\begin{tikzpicture}
\begin{axis}[
  cqr/timings,
  ymax=100,ytick={0,20,...,100},
  title={Tolerance: $\kappa_2(Q_{\mathrm{p}})<10^{8}$}
]
  \AddTimingSeries{\cqrbardata}
\end{axis}
\end{tikzpicture}
\endgroup
\caption{Total runtime required by \cref{alg:iterprec} to achieve
$\cond_2(Q_{\mathrm{p}})<10^8$ for a matrix $X$ with
$\cond_2(X)=10^{15}$.
The different sections of each bar represent the individual runtimes of different components.
The total number of iterations is shown above the bar.
The cost of estimating $\cond_2(Q_{\mathrm{p}})$ is excluded.}
\label{fig:compbar1e8}
\end{figure}

Next, we compare the runtime of the ten preconditioning algorithms on the test
matrices in~\eqref{eq:test-matrix} for $n = 100$ and $\kappa=10^{15}$.
\Cref{fig:compbar1e1,fig:compbar1e8} report the time spent by
\cref{alg:iterprec}, broken down by operation, for $\tau=10$ and $\tau=10^8$,
respectively.

For $\tau = 10$,
RCQR with a sparse sign sketch is the fastest preconditioner, followed by
mpLUCQR, which is the fastest among the three proposed alternatives.
For $\tau = 10^8$, mpLUCQR becomes the fastest preconditioner,
while mpEVDQR becomes comparable to RCQR with a
sparse sign sketch.

We now consider the timing breakdown. On the computing platform used in our
experiments, the peak performance of matrix--matrix multiplication in binary32
arithmetic is approximately $64\times$ that of binary64 arithmetic. This speedup
is most significant in the computations of the Gram matrices $X\trans X$ and
$L\trans L$ and of the sketched matrix. Comparing SCQR with mpSCQR and mpEVDQR
shows that it is approximately $22\times$ faster to form $X\trans X$ in low
precision compared with working precision. A similar comparison between LUCQR
and mpLUCQR shows an approximately $45\times$ acceleration in forming
$L\trans L$ in binary32 arithmetic. Similarly, comparing RCQR and mpRCQR with
Gaussian sketching shows that forming $\Theta X$ in low precision is
approximately $25\times$ faster than in working precision.

Comparing LUCQR with mpLUCQR shows that a binary32 LU factorization is approximately $4.6\times$ faster than a
binary64 one and $5.6\times$ faster
than binary64 matrix--matrix multiplication.
Similarly, the binary32 Householder QR factorization in mpHQR is $2.8\times$
faster than a binary64 matrix multiplication.
The only exception is RCQR with sparse sign sketching, since the
sketching matrix is sparse and the main bottleneck in forming the sketched
matrix is likely to be data movement.

We conclude by noting that a speedup in constructing the preconditioner does not
guarantee a reduction in the overall runtime, because each iteration also
requires a binary64 update of the preconditioned matrix for~\eqref{eq:XQRerror}
to hold. Therefore, a cheaper preconditioner may yield worse performance overall
if more iterations are required. In \cref{fig:compbar1e8}, for example, mpSCQR
is slower than SCQR because it requires three iterations rather than one. Among
the three proposed methods, mpSCQR and mpEVDQR have relatively low cost per
iteration, but they only reduce the condition number by a modest factor of
$O(\sqrt{\ul})$. Conversely, a preconditioner that yields large reduction in
condition number does not necessarily lead to a shorter overall runtime. For
example, a single application of LUCQR reduces the condition number from
$10^{15}$ to below $10$, but factorizing $X$ and forming $L\trans L$ in binary64
arithmetic are costly operations that make it slower than mpLUCQR for both
values of $\tau$ considered.

\section{Conclusions and future work}

We proposed three mixed-precision preconditioners for CholeskyQR.
These are based on shifted Cholesky factorization, eigendecomposition, and partially pivoted LU factorization.
We bounded the rate at which these preconditioners reduce the condition number and we
compared them against several existing alternatives using GPU implementations.

Our analysis gives a new interpretation of CholeskyQR preconditioning as an iterative refinement scheme.
Existing analyses~\cite{fkny20,too20} typically show that, under certain assumptions,
one preconditioning step reduces the condition number sufficiently for CholeskyQR or CholeskyQR2 to succeed.
Our results also quantify the reduction achieved by our preconditioners, which is approximately
$O(\sqrt{\ul})$ for mpSCQR and mpEVDQR and $O(\ul)$ for mpLUCQR. However, we argue that when one
preconditioning step is not sufficient, repeated applications will progressively reduce the
condition number until a desired condition number is reached.

The numerical behavior of the proposed preconditioners agrees with our theoretical results, and
our experiments show that they can be considerably more efficient than the mono-precision
shiftedCholeskyQR and mono-precision LU--CholeskyQR preconditioners, achieving performance on par with that of randomized preconditioners.

Several questions still remain.
First, the convergence result for mpLUCQR relies on the relatively strong
assumption~\eqref{eq:mplucqr-assump-chi}, and it would be valuable to either
weaken this assumption or modify the algorithm so that a stronger result can be obtained.
Second, estimating $\cond_2(Q_{\mathrm{p}})$ at each iteration of \cref{alg:iterprec} to
decide when to terminate the algorithms can be very expensive.
It may be possible to replace this computation with a cheaper alternative, such as
$\cond_2(R_{\mathrm{p}})$, for example, but further analysis is required to determine
when this cheaper alternative is sufficiently reliable.

\section*{Acknowledgments}
This work used the HPC Hardware Lab which is part of the COSMA facility which includes the DiRAC@Durham facility managed by the Institute for Computational Cosmology on behalf of the STFC DiRAC HPC Facility (www.dirac.ac.uk). The equipment was funded by BEIS capital funding via STFC capital grants ST/K00042X/1, ST/P002293/1, ST/R002371/1 and ST/S002502/1, Durham University, and STFC operations grant ST/R000832/1. DiRAC is part of the National e-Infrastructure. The HPC Hardware Lab is funded from various sources including the UK ExCALIBUR programme.
GPT-6 was used to assist in deriving the bounds in~\cref{thm:mpscqr-convergence,thm:mplucqr-convergence} and to refine the code used for the experiments.

\appendix

\section{Algorithms for comparison}
\label{sec:otheralg}
In this section, we introduce the two preconditioners for CholeskyQR that we consider in \cref{sec:experiments}.

\subsection{Mixed-precision Householder QR}

\begin{algorithm2e}[t]
\caption{Mixed-precision Householder QR (mpHQR)}
\label{alg:MPHQR}
\Indm
\KwIn{$X\in\mathbb{R}^{m\times n}$; unit roundoffs $u\le \ul$.}
\KwOut{CholeskyQR preconditioner $R_{\mathrm{p}}$;
  preconditioned matrix $Q_{\mathrm{p}}=XR_{\mathrm{p}}^{-1}$.}
\Indp
\Fn{$[Q_{\mathrm{p}},R_{\mathrm{p}}]=\texttt{MPHQR}(X)$}{
    Compute a Householder QR factorization $Q_{\mathrm{p}}R_{\mathrm{p}} = X$ in precision $\ul$.\;
    Convert $R_{\mathrm{p}}$ to precision $u$.\;
    Compute $Q_{\mathrm{p}}\gets XR_{\mathrm{p}}^{-1}$
    in precision $u$.\label{alg-step:hqr_update}\;
}
\end{algorithm2e}

\emph{Mixed-precision Householder QR} (mpHQR) in \cref{alg:MPHQR} was first proposed by Kressner, Ma, and Shao~\cite{kms23} for orthogonalizing a basis of the search space of LOBPCG.
It first computes a Householder QR factorization in low precision and
retains only its upper triangular factor.
It then forms the preconditioned matrix through a working-precision triangular
solve.

On \cref{alg-step:hqr_update}, $Q_{\mathrm{p}}$ is computed as $XR_{\mathrm{p}}^{-1}$ instead of the product
of Householder reflectors. This is necessary because, in general, the product of low-precision
Householder reflectors
cannot satisfy ~\eqref{eq:XQRerror} even if stored in high precision.

\subsection{Randomized CholeskyQR}

\emph{Randomized CholeskyQR} (RCQR)~\cite{bala22,hsby25} is a preconditioner for CholeskyQR based on sketching.
Let $\Theta\in\R^{d\times m}$, with $d\geq n$, be an
$\varepsilon$-subspace embedding for $\operatorname{range}(X)$, that is, such that
\begin{equation}
(1-\varepsilon)\norm{Xz}_2^2
\leq \norm{\Theta Xz}_2^2
\leq (1+\varepsilon)\norm{Xz}_2^2,
\qquad z\in\R^n.
\label{eq:subspace_embedding}
\end{equation}
RCQR forms the sketched matrix \mbox{$X_{\Theta}=\Theta X$}, computes a Householder
QR factorization $X_{\Theta}=Q_{\Theta}R_{\mathrm{p}}$, and sets
$Q_{\mathrm{p}}=XR_{\mathrm{p}}^{-1}$.  Since $\Theta Q_{\mathrm{p}}=Q_{\Theta}$,
condition \eqref{eq:subspace_embedding} gives
\begin{equation*}
\cond_2(Q_{\mathrm{p}})\leq
\sqrt{\frac{1+\varepsilon}{1-\varepsilon}}.
\end{equation*}
Thus, with a suitable sketching matrix $\Theta$, $Q_{\mathrm{p}}$ is well
conditioned with high probability,
and a subsequent CholeskyQR step can then recover orthogonality in working precision.
\footnote{In fact, the finite-precision analyses of $\cond_2(\wh Q_{\mathrm{p}})$ are also given in
\cite{bala22,hsby25}.
However, randomized algorithms are not the main focus of our work and are therefore not discussed in detail here.}
Here, we consider two choices for $\Theta$: a Gaussian matrix~\cite{matr20} and a sparse sign matrix~\cite{clwo17}.

A variant of this preconditioner that performs the sketching in low precision, mpRCQR in \cref{alg:MPRCQR}, has recently been proposed~\cite{gci26}.
Setting $\ul=u$ in \cref{alg:MPRCQR} recovers RCQR.

\begin{algorithm2e}[t]
\caption{Mixed-precision randomized CholeskyQR (mpRCQR)}
\label{alg:MPRCQR}
\Indm
\KwIn{$X\in\mathbb{R}^{m\times n}$;
sketching matrix $\Theta\in\mathbb{R}^{d\times m}$;
unit roundoffs $u\le \ul$.}
\KwOut{CholeskyQR preconditioner $R_{\mathrm{p}}$;
preconditioned matrix $Q_{\mathrm{p}}=XR_{\mathrm{p}}^{-1}$.}
\Indp
\Fn{$[Q_{\mathrm{p}},R_{\mathrm{p}}]=\texttt{MPRCQR}(X)$}{
  Compute the sketch $X_{\Theta}=\Theta X$ in precision $\ul$.\;
  Compute a Householder QR factorization $Q_{\mathrm{\Theta}}R_{\mathrm{p}} = X_{\Theta}$ in precision $u$.\;
  Compute $Q_{\mathrm{p}}\gets XR_{\mathrm{p}}^{-1}$ in precision $u$.\;
}
\end{algorithm2e}

\section{Proof of Theorem~\ref{thm:mpscqr-convergence}}
\label{sec:conv_mpscqr_proof}
To simplify the notation, throughout this section we write $Q=Q_{\mathrm{p}}$ and $R=R_{\mathrm{p}}$.
Push c
Using the error bounds of~\cite{jeru13,ruje14}, the componentwise error bounds for \cref{alg:MPSCQR} are
\begin{align}
\wh G-X\trans X&=E_1,
&\abs{E_1} &\le (m+3)\ul\abs{X\trans} \abs{X},\label{eq:mpscqr-E1}\\
\wh G_s-(\wh G+sI_n)&=E_2,
&\abs{E_2} &\le u\abs{\wh G+sI_n},\label{eq:mpscqr-E2}\\
\wh R\trans\wh R-\wh G_s&=E_3,
&\abs{E_3} &\le (n+1)u\abs{\wh R\trans} \abs{\wh R},\label{eq:mpscqr-E3}\\
\wh Q\wh R-X&=E_4,
&\abs{E_4} &\le nu\abs{\wh Q} \abs{\wh R},\label{eq:mpscqr-E4}
\end{align}
where $\wh G_s$ is defined as the computed $\wh G+sI_n$.
The bound~\eqref{eq:mpscqr-E1}, which includes the error in converting $X$ to the low-precision matrix $X_\ell$,
can be obtained from
\begin{align*}
\abs{E_1}&\le m\ul\abs{X_\ell\trans}\abs{X_\ell}
+\abs{X_\ell\trans X_\ell-X\trans X}\\
&\le\bigl(m\ul(1+\ul)^2+2\ul+\ul^2\bigr)\abs{X\trans}\abs{X}\\
&\le\bigl((m+3)\ul\bigr)\abs{X\trans}\abs{X}
\end{align*}
by using~\cite[Thm.~4.2]{jeru13} and~\eqref{eq:mpscqr-small-dim} in the second and third inequality, respectively. The bounds~\eqref{eq:mpscqr-E2} and \eqref{eq:mpscqr-E4} follow from~\cite[eq.~(2.4)]{high02} and~\cite[Thm.~4.2]{jeru13}, respectively, while \eqref{eq:mpscqr-E3} can be found in~\cite[Thm.~4.4]{ruje14}.

From~\eqref{eq:mpscqr-E1}--\eqref{eq:mpscqr-E4}, we can readily derive the normwise error bounds
\begin{align}
\norm{E_1}_2&\le (m+3)n\ul\norm{X}_2^2,\label{eq:mpscqr-E1-norm}\\
\norm{E_2}_2&\le u\bigl(\norm{\wh G}_2+s\bigr),\label{eq:mpscqr-E2-norm}\\
\norm{E_3}_2&\le n(n+1)u\norm{\wh R}_2^2,\label{eq:mpscqr-E3-norm}\\
\norm{E_4}_2&\le n^2u\norm{\wh Q}_2\norm{\wh R}_2,\label{eq:mpscqr-E4-norm}
\end{align}
where in~\eqref{eq:mpscqr-E2-norm} we have used the fact that $E_2$ is diagonal.
We further assume that $E_1$, $E_2$, and $E_3$ are symmetric, as they are produced by symmetric operations.

We begin by estimating the error in the Cholesky factorization.
\begin{lemma}[Error in the Cholesky factorization]
\label{lem:mpscqr-small-perturbation}
If~\eqref{eq:mpscqr-shift-range} is satisfied, then the Cholesky factorization will not break down and
\begin{equation}
\norm{E_1+E_2+E_3}_2\le 0.1s.\label{eq:mpscqr-total-error}
\end{equation}
Consequently,
\begin{equation}
\label{eq:mpscqr-loewner}
X\trans X+0.9sI_n
\preceq \wh R\trans\wh R
\preceq X\trans X+1.1sI_n,
\end{equation}
where $\preceq$ denotes the Loewner partial order.
\end{lemma}

\begin{proof}
We first verify that the shift prevents the Cholesky factorization from breaking down.
By~\eqref{eq:mpscqr-E1-norm}, \eqref{eq:mpscqr-E2-norm},
\eqref{eq:mpscqr-small-dim}, and $s/\norm{X}_2^2\le1/100$, we have
\begin{equation}\label{eq:mpscqr-E2-final}
\norm{E_2}_2
\le u\bigl(\norm{\wh G}_2+s\bigr)
\le u\bigl((1+(m+3)n\ul)\norm{X}_2^2+s\bigr)
\le1.1u\norm{X}_2^2.
\end{equation}
Let $G_0=\wh G+E_2+\bigl((m+3)n\ul+1.1u\bigr)\norm{X}_2^2I_n$.
By~\eqref{eq:mpscqr-E1-norm} and the bound on $E_2$, we have
\begin{equation*}
G_0=E_1+E_2+X\trans X+
\bigl((m+3)n\ul+1.1u\bigr)\norm{X}_2^2I_n\succeq X\trans X,
\end{equation*}
so $G_0$ is positive semidefinite.
Moreover, \eqref{eq:mpscqr-E1}~gives
\begin{equation*}
\trace(E_1)
\le \trace(\abs{E_1})
\le (m+3)\ul\trace(\abs{X\trans} \abs{X})
= (m+3)\ul\norm{X}_F^2
\le (m+3)n\ul\norm{X}_2^2,
\end{equation*}
and therefore
\begin{equation}\label{eq:trg0-bound}
\begin{aligned}
\trace(G_0)
&=\trace(E_1)+\trace(E_2)+\trace(X\trans X)
+\bigl((m+3)n\ul+1.1u\bigr)\norm{X}_2^2\cdot\trace(I_n)\\
&\le\bigl(n+(m+3)n(n+1)\ul+2.2nu\bigr)\norm{X}_2^2\\
&\le1.1n\norm{X}_2^2,
\end{aligned}
\end{equation}
where the last inequality follows from~\eqref{eq:mpscqr-small-dim}.
The matrix to which Cholesky factorization is applied is
\begin{equation*}
\wh G_s
=G_0+\left(s-\bigl((m+3)n\ul+1.1u\bigr)\norm{X}_2^2\right)I_n.
\end{equation*}
Thus, by~\cite[sect.~4.1]{fkny20}, the Cholesky factorization will not break down if
\begin{equation*}
s\ge\bigl((m+3)n\ul+1.1u\bigr)\norm{X}_2^2
+2.2(n+1)u\cdot\trace(G_0).
\end{equation*}
Using~\eqref{eq:trg0-bound}, we obtain
\begin{equation*}
\begin{aligned}
&\bigl((m+3)n\ul+1.1u\bigr)\norm{X}_2^2
+2.2(n+1)u\cdot\trace(G_0)\\
&\qquad\le\bigl((m+3)n\ul+1.1u+2.42n(n+1)u\bigr)\norm{X}_2^2\\
&\qquad\le11\bigl((m+3)n\ul+n(n+3)u\bigr)\norm{X}_2^2\le s,
\end{aligned}
\end{equation*}
which shows that $s$ in~\eqref{eq:mpscqr-shift-range} guarantees completion of the Cholesky factorization.

Next, we estimate $\norm{E_1+E_2+E_3}_2$.
Using $\wh R\trans \wh R=\wh G + sI_n + E_2+E_3$,
we obtain
\begin{equation}
  \label{eq:rhat-bound}
\norm{\wh R}_2^2
\le \norm{X}_2^2+s+\norm{E_1}_2+\norm{E_2}_2+\norm{E_3}_2.
\end{equation}
Substituting~\eqref{eq:rhat-bound} into~\eqref{eq:mpscqr-E3-norm} and using
\eqref{eq:mpscqr-E1-norm} and \eqref{eq:mpscqr-E2-final} gives
\begin{equation}
\label{eq:mpscqr-E3-final}
\begin{aligned}
\norm{E_3}_2
&\le
\frac{n(n+1)u}{1-n(n+1)u}
\left(1+\frac{s}{\norm{X}_2^2}+(m+3)n\ul+1.1u\right)
\norm{X}_2^2\\
&\le 1.1n(n+1)u\norm{X}_2^2,
\end{aligned}
\end{equation}
where the last inequality follows from~\eqref{eq:mpscqr-small-dim} and $s/\norm{X}_2^2\le1/100$.
Combining~\eqref{eq:mpscqr-E1-norm}, \eqref{eq:mpscqr-E2-final}, and \eqref{eq:mpscqr-E3-final}, we obtain
\begin{equation*}
\begin{aligned}
 \norm{E_1+E_2+E_3}_2
 &\le\bigl((m+3)n\ul+1.1u+1.1n(n+1)u\bigr)\norm{X}_2^2\\
 &\le1.1\bigl((m+3)n\ul+n(n+3)u\bigr)\norm{X}_2^2
 \le0.1s.
\end{aligned}
\end{equation*}
Since $E_1+E_2+E_3$ is symmetric, \eqref{eq:mpscqr-loewner} follows.
\end{proof}

If we assume that the solution of the triangular system doe not introduce any additional error,
we can use \cref{lem:mpscqr-small-perturbation} to estimate $\cond_2(Q)$ .
\begin{lemma}[The ideally preconditioned matrix]
\label{lem:mpscqr-exact-solve}
Let $Y=X\wh R^{-1}$.
If~\eqref{eq:mpscqr-shift-range} is satisfied, then
\begin{align}
\cond_2(Y)
&\le\sqrt{\frac{1+1.1s\cond_2(X)^2/\norm{X}_2^2}{1+0.9s/\norm{X}_2^2}},\label{eq:mpscqr-Ycond}\\
\cond_2(\wh R)
&\le\sqrt{\frac{1+1.1s/\norm{X}_2^2}{0.9s/\norm{X}_2^2}}.\label{eq:mpscqr-Rcond}
\end{align}
\end{lemma}

\begin{proof}
For every nonzero $z\in\R^n$, \eqref{eq:mpscqr-loewner}~gives
\begin{equation*}
\norm{Xz}_2^2+0.9s\norm{z}_2^2
\le\norm{\wh Rz}_2^2
\le\norm{Xz}_2^2+1.1s\norm{z}_2^2,
\end{equation*}
or equivalently
\begin{equation}
\frac{\norm{Xz}_2^2}{\norm{Xz}_2^2+1.1s\norm{z}_2^2}
\le\frac{\norm{Xz}_2^2}{\norm{\wh Rz}_2^2}
\le\frac{\norm{Xz}_2^2}{\norm{Xz}_2^2+0.9s\norm{z}_2^2}.\label{eq:mpscqr-Rayleigh}
\end{equation}
Since $\wh R$ is nonsingular, we have
\begin{equation}\label{eq:min-max}
\max_{z\neq0}\frac{\norm{Xz}_2^2}{\norm{\wh Rz}_2^2}
=\max_{y\neq0}\frac{\norm{Yy}_2^2}{\norm{y}_2^2}
=\sigma_{\max}(Y)^2,\qquad
\min_{z\neq0}\frac{\norm{Xz}_2^2}{\norm{\wh Rz}_2^2}
=\min_{y\neq0}\frac{\norm{Yy}_2^2}{\norm{y}_2^2}
=\sigma_{\min}(Y)^2.
\end{equation}
Dividing the numerator and denominator of the lower and upper bounds
in~\eqref{eq:mpscqr-Rayleigh} by $\norm{z}_2^2$ and using~\eqref{eq:min-max} yields
\begin{equation*}
\frac{\sigma_{\min}(X)^2}{\sigma_{\min}(X)^2+1.1s}
\le\sigma_{\min}(Y)^2
\le\sigma_{\max}(Y)^2
\le\frac{\sigma_{\max}(X)^2}{\sigma_{\max}(X)^2+0.9s},
\end{equation*}
which proves~\eqref{eq:mpscqr-Ycond}.

Finally, it also follows from~\eqref{eq:mpscqr-loewner} that
\begin{equation*}
\sigma_{\max}(\wh R)^2\le\norm{X}_2^2+1.1s,\qquad
\sigma_{\min}(\wh R)^2\ge0.9s,
\end{equation*}
which proves \eqref{eq:mpscqr-Rcond}.
\end{proof}

With the help of \cref{lem:mpscqr-small-perturbation,lem:mpscqr-exact-solve},
we can finally prove \cref{thm:mpscqr-convergence}.
\begin{proof}[Proof of \cref{thm:mpscqr-convergence}]
Let
\begin{equation*}
\alpha=\frac{s}{\norm{X}_2^2}.
\end{equation*}
By~\eqref{eq:mpscqr-shift-range} and~\eqref{eq:mpscqr-Rcond}, we have
\begin{equation*}
n^2u\cond_2(\wh R)
\le
\frac{\sqrt{\alpha(1+1.1\alpha)}}{11\sqrt{0.9}}
<0.01,
\end{equation*}
where we used
\begin{equation}
  \label{eq:bound-alpha}
  11n^2u\le\alpha\le\frac{1}{100},
\end{equation}
which is implied by~\eqref{eq:mpscqr-shift-range}.
Hence $1-n^2u\cond_2(\wh R)>0$.

Multiplying $\wh Q\wh R-X=E_4$ by $\wh R^{-1}$ on the right gives
\begin{equation*}
\wh Q-Y=E_4\wh R^{-1}.
\end{equation*}
By~\eqref{eq:mpscqr-E4-norm}, we have
\begin{equation*}
\norm{\wh Q}_2\le\norm{Y}_2+\norm{E_4\wh R^{-1}}_2
\le \norm{Y}_2+n^2u\cond_2(\wh R)\norm{\wh Q}_2.
\end{equation*}
Since $1-n^2u\cond_2(\wh R)>0$, this leads to
\begin{equation}
\label{eq:mpscqr-sigmamaxQp}
\norm{\wh Q}_2\le\frac{\norm{Y}_2}{1-n^2u\cond_2(\wh R)}.
\end{equation}
Weyl's theorem then yields
\begin{equation}
\label{eq:mpscqr-sigmaminQp}
\begin{aligned}
\sigma_{\min}(\wh Q)&\ge\sigma_{\min}(Y)-\norm{E_4\wh R^{-1}}_2\\
&\ge\sigma_{\min}(Y)-n^2u\cond_2(\wh R)\norm{\wh Q}_2\\
&\ge\sigma_{\min}(Y)-\frac{n^2u\cond_2(\wh R)}
{1-n^2u\cond_2(\wh R)}\norm{Y}_2.
\end{aligned}
\end{equation}
Combining~\eqref{eq:mpscqr-sigmamaxQp} and~\eqref{eq:mpscqr-sigmaminQp}, we obtain
\begin{equation}
\label{eq:mpscqr-transfer}
\cond_2(\wh Q)\le\frac{\cond_2(Y)}
{1-n^2u\cond_2(\wh R)\bigl(1+\cond_2(Y)\bigr)}.
\end{equation}

It remains to bound the denominator of~\eqref{eq:mpscqr-transfer}.
Using \cref{lem:mpscqr-exact-solve}, \eqref{eq:bound-alpha}, and $\sqrt{1+t^2}\le1+t$ for $t\ge0$, we obtain
\begin{equation*}
\begin{aligned}
n^2u\cond_2(\wh R)\bigl(1+\cond_2(Y)\bigr)&\le
n^2u\sqrt{\frac{1+1.1\alpha}{0.9\alpha}}
\left(1+\sqrt{1+1.1\alpha\cond_2(X)^2}\right)\\
&<
2.12\frac{n^2u}{\sqrt{\alpha}}
+1.12n^2u\cond_2(X)\\
&<\frac{2.12}{110}+\frac{1.12}{6}\\
&<0.21,
\end{aligned}
\end{equation*}
where in the last inequality we used~\eqref{eq:mpscqr-assump-cond}.
Hence the denominator in \eqref{eq:mpscqr-transfer} is positive, and substituting \eqref{eq:mpscqr-Ycond} into \eqref{eq:mpscqr-transfer} proves \eqref{eq:mpscqr-main-bound}.
\end{proof}

\section{Proof of Theorem~\ref{thm:mplucqr-convergence}}
\label{sec:conv_mplucqr_proof}

To simplify the notation, throughout this section, we write $Q=Q_{\mathrm{p}}$ and $R=R_{\mathrm{p}}$.

The componentwise error bounds of \cref{alg:MPLUCQR} are
\begin{equation*}
\begin{aligned}
\wh G-\wh L\trans\wh L&=E_2,
\qquad\abs{E_2}\le m\ul\abs{\wh L\trans}\abs{\wh L},\\
\wh S\trans\wh S-\wh G&=E_3,
\qquad\abs{E_3}\le(n+1)u\abs{\wh S\trans}\abs{\wh S},\\
\wh R-\wh S\wh U&=E_4,
\qquad\abs{E_4}\le nu\abs{\wh S}\abs{\wh U},\\
\wh Q\wh R-X&=E_5,
\qquad\abs{E_5}\le nu\abs{\wh Q}\abs{\wh R},
\end{aligned}
\end{equation*}
and the corresponding normwise bounds are
\begin{align}
\norm{E_2}_2&\le mn\ul\norm{\wh L}_2^2,
\label{eq:mplucqr-E2-norm}\\
\norm{E_3}_2&\le n(n+1)u\norm{\wh S}_2^2,
\label{eq:mplucqr-E3-norm}\\
\norm{E_4}_2&\le n^2u\norm{\wh S}_2\norm{\wh U}_2,
\label{eq:mplucqr-E4-norm}\\
\norm{E_5}_2&\le n^2u\norm{\wh Q}_2\norm{\wh R}_2.
\label{eq:mplucqr-E5-norm}
\end{align}
There is no $E_1$ since the error of the first step, including the rounding error of converting
$X$ to low precision and LU factorization, is represented by $\chi_{\min}$ and
$\chi_{\max}$ in \cref{thm:mplucqr-convergence}; see also \cref{rmk:mplucqr-assump-chi-1}.

Using assumptions~\eqref{eq:mplucqr-assump-u-L} and~\eqref{eq:mplucqr-assump-ul}, one can prove that the Cholesky
factorization of $\wh G$ will not break down.
We do not repeat the proof, which is analogous to that of \cref{lem:mpscqr-small-perturbation}.

\begin{lemma}[Extreme singular values of $\wh S$ and $\wh R$]
\label{lem:mplucqr-factor-bounds}
Under assumptions~\eqref{eq:mplucqr-assump-u-L}, \eqref{eq:mplucqr-assump-u-LU},
and~\eqref{eq:mplucqr-assump-ul}, we have
\begin{align}
\norm{\wh S}_2
&\le1.02\norm{\wh L}_2,\label{eq:mplucqr-Snorm}\\
\norm{\wh S^{-1}}_2
&\le1.02\sigma_{\min}(\wh L)^{-1},\label{eq:mplucqr-Sinvnorm}\\
\norm{\wh R^{-1}}_2
&\le1.04\sigma_{\min}(\wh L)^{-1}\sigma_{\min}(\wh U)^{-1},\label{eq:mplucqr-Rinvnorm}\\
\cond_2(\wh R)
&\le1.066\cond_2(\wh L)\cond_2(\wh U).\label{eq:mplucqr-Rcond}
\end{align}
\end{lemma}

\begin{proof}
We first prove the bounds for $\wh S$.
Assumptions~\eqref{eq:mplucqr-assump-u-L} and~\eqref{eq:mplucqr-assump-ul} imply that
$mn\ul\le1/64$ and $n(n+1)u\le1/64$.
Therefore, substituting~\eqref{eq:mplucqr-E2-norm} and
\eqref{eq:mplucqr-E3-norm} into $\wh S\trans\wh S=\wh L\trans\wh L+E_2+E_3$ gives
\begin{equation*}
\norm{\wh S}_2^2
\le
\frac{1+mn\ul}{1-n(n+1)u}\norm{\wh L}_2^2\le\frac{ 65}{63}\norm{\wh L}_2^2.
\end{equation*}
On the other hand, by Weyl's theorem, the smallest singular value of $\wh S$ satisfies
\begin{equation*}
\sigma_{\min}(\wh S)^2
\ge\sigma_{\min}(\wh L)^2-\norm{E_2+E_3}_2
\ge\frac{61}{63}\sigma_{\min}(\wh L)^2,
\end{equation*}
where, by~\eqref{eq:mplucqr-assump-u-L} and~\eqref{eq:mplucqr-assump-ul}, the last inequality comes from
\begin{equation}
\label{eq:mplucqr-E23-bound}
\begin{aligned}
\norm{E_2+E_3}_2
&\le\norm{E_2}_2+\norm{E_3}_2\\
&\le\left(mn\ul+\frac{65}{63}n(n+1)u\right)
\norm{\wh L}_2^2\\
&\le\left(mn\ul\cond_2^2(\wh L)+\frac{65}{63}n(n+1)u\cond_2^2(\wh L)\right)
\sigma_{\min}(\wh L)^2\\
&\le\frac{2}{63}\sigma_{\min}(\wh L)^2.
\end{aligned}
\end{equation}

For $\wh R$, if we apply Weyl's theorem together with ~\eqref{eq:mplucqr-E4-norm}, \eqref{eq:mplucqr-Snorm}, and~\eqref{eq:mplucqr-Sinvnorm}, we obtain
\begin{equation*}
\begin{aligned}
\sigma_{\min}(\wh R)
&\ge\sigma_{\min}(\wh S)\sigma_{\min}(\wh U)-n^2u\norm{\wh S}_2\norm{\wh U}_2\\
&\ge\bigg(\sqrt{\frac{61}{63}}
-\sqrt{\frac{65}{63}}n^2u\cond_2(\wh L)\cond_2(\wh U)\bigg)
\sigma_{\min}(\wh L)\sigma_{\min}(\wh U)\\
&\ge\bigg(\sqrt{\frac{61}{63}}
-\frac{1}{64}\sqrt{\frac{65}{63}}\bigg)
\sigma_{\min}(\wh L)\sigma_{\min}(\wh U),
\end{aligned}
\end{equation*}
where, in the last inequality, we used~\eqref{eq:mplucqr-assump-u-LU}
to bound $n^2u\cond_2(\wh L)\cond_2(\wh U)$.
Finally, dividing
\begin{equation*}
\norm{\wh R}_2\le(1+n^2u)\norm{\wh S}_2\norm{\wh U}_2
\le\frac{65}{64}\cdot\sqrt{\frac{65}{63}}\norm{\wh L}_2\norm{\wh U}_2
\end{equation*}
by $\sigma_{\min}(\wh R)$ yields~\eqref{eq:mplucqr-Rcond}.
\end{proof}

\begin{lemma}[Near-isometry of the LU-based preconditioner]
\label{lem:mplucqr-LUR-condition}
Under assumptions~\eqref{eq:mplucqr-assump-u-L},
\eqref{eq:mplucqr-assump-u-LU}, and~\eqref{eq:mplucqr-assump-ul}, we have
\begin{equation}
\label{eq:mplucqr-LUR-bounds}
0.965
\le\sigma_{\min}(\wh L\wh U\wh R^{-1})
\le\sigma_{\max}(\wh L\wh U\wh R^{-1})
\le1.034.
\end{equation}
\end{lemma}

\begin{proof}
We first estimate the perturbation introduced when
$\wh S\wh U$ is formed in working precision.
By~\eqref{eq:mplucqr-E4-norm},
\eqref{eq:mplucqr-Snorm}, and~\eqref{eq:mplucqr-Rinvnorm}, we have
\begin{equation*}
\begin{aligned}
\norm{E_4\wh R^{-1}}_2
\le n^2u\norm{\wh S}_2\norm{\wh U}_2\norm{\wh R^{-1}}_2
\le1.02\cdot1.04 n^2u\cond_2(\wh L)\cond_2(\wh U).
\end{aligned}
\end{equation*}
Assumption~\eqref{eq:mplucqr-assump-u-LU} therefore gives
\begin{equation}
\label{eq:mplucqr-E4Rinv}
\norm{E_4\wh R^{-1}}_2
\le\frac{1.02\cdot1.04}{64}
=0.016575
<0.0166.
\end{equation}

Next, multiplying
$\wh R=\wh S\wh U+E_4$ on the right by $\wh R^{-1}$
gives
\begin{equation*}
\wh S\wh U\wh R^{-1}
=I_n-E_4\wh R^{-1}.
\end{equation*}
Since $\wh S$ is nonsingular by \cref{lem:mplucqr-factor-bounds}, it follows that
\begin{equation*}
\wh U\wh R^{-1}
=\wh S^{-1}\bigl(I_n-E_4\wh R^{-1}\bigr).
\end{equation*}
Using~\eqref{eq:mplucqr-Sinvnorm} and
\eqref{eq:mplucqr-E4Rinv}, we obtain
\begin{equation*}\norm{\wh U\wh R^{-1}}_2\le
1.02(1+0.0166)\sigma_{\min}(\wh L)^{-1}.
\end{equation*}
Then, using~\eqref{eq:mplucqr-E23-bound}, we obtain
\begin{equation}
\label{eq:mplucqr-E23UR-bound}
\begin{aligned}
\norm{\wh U\wh R^{-1}}_2^2\bigl(\norm{E_2}_2+\norm{E_3}_2\bigr)&\le1.02^2(1+0.0166)^2\sigma_{\min}(\wh L)^{-2}
\cdot\frac{2}{63}\sigma_{\min}(\wh L)^2\\
&=\frac{2\cdot1.02^2}{63}(1+0.0166)^2\\
&<0.0342.
\end{aligned}
\end{equation}

We now estimate the deviation from orthogonality of $\wh L\wh U\wh R^{-1}$.
Note that
\begin{equation*}
\ifsimax
\begin{aligned}
(\wh L\wh U\wh R^{-1})\trans
(\wh L\wh U\wh R^{-1})-I_n=
&~-\wh R\transinv E_4\trans
-E_4\wh R^{-1}
+\wh R\transinv E_4\trans E_4\wh R^{-1}\\
&~-\wh R\transinv\wh U\trans(E_2+E_3)
 \wh U\wh R^{-1}.
\end{aligned}
\else
(\wh L\wh U\wh R^{-1})\trans
 (\wh L\wh U\wh R^{-1})-I_n =
-\wh R\transinv E_4\trans
-E_4\wh R^{-1}
+\wh R\transinv E_4\trans E_4\wh R^{-1}
-\wh R\transinv\wh U\trans(E_2+E_3)
 \wh U\wh R^{-1}.
\fi
\end{equation*}
Using~\eqref{eq:mplucqr-E4Rinv} and~\eqref{eq:mplucqr-E23UR-bound}, we obtain
\begin{equation*}
\ifsimax
\begin{aligned}
&~\norm{
(\wh L\wh U\wh R^{-1})\trans
(\wh L\wh U\wh R^{-1})-I_n}_2\\
\le&~2\norm{E_4\wh R^{-1}}_2
+\norm{E_4\wh R^{-1}}_2^2
+\norm{\wh U\wh R^{-1}}_2^2
\bigl(\norm{E_2}_2+\norm{E_3}_2\bigr)\\
<&~0.0677.
\end{aligned}
\else
\norm{
(\wh L\wh U\wh R^{-1})\trans
(\wh L\wh U\wh R^{-1})-I_n}_2
\le2\norm{E_4\wh R^{-1}}_2
+\norm{E_4\wh R^{-1}}_2^2
+\norm{\wh U\wh R^{-1}}_2^2
 \bigl(\norm{E_2}_2+\norm{E_3}_2\bigr)<0.0677.
\fi
\end{equation*}
Therefore,
\begin{equation*}
\sigma_{\min}(\wh L\wh U\wh R^{-1})
\ge\sqrt{0.9323}>0.965,
\end{equation*}
and
\begin{equation*}
\sigma_{\max}(\wh L\wh U\wh R^{-1})
\le\sqrt{1.0677}<1.034,
\end{equation*}
which proves~\eqref{eq:mplucqr-LUR-bounds}.
\end{proof}

\begin{lemma}[Ideally preconditioned matrix]
\label{lem:mplucqr-XR-condition}
Let $Y=X\wh R^{-1}$.  Under assumptions
\eqref{eq:mplucqr-assump-u-L},
\eqref{eq:mplucqr-assump-u-LU}, and~\eqref{eq:mplucqr-assump-ul},
\begin{equation}
\label{eq:mplucqr-XR-bounds}
0.965\cdot\chi_{\min}
\le\sigma_{\min}(Y)
\le\sigma_{\max}(Y)
\le1.034\cdot\chi_{\max},
\end{equation}
and consequently
\begin{equation*}
\cond_2(Y)
\le1.072\cdot\frac{\chi_{\max}}{\chi_{\min}}.
\end{equation*}
\end{lemma}

\begin{proof}
By the definition of $\chi_{\min}$ and $\chi_{\max}$ in \eqref{eq:chimin-chimax-def},
for every $z\in\mathbb{R}^n$, we have
\begin{equation}
\label{eq:mplucqr-chi-vector-bound}
\chi_{\min}\norm{\wh L\wh Uz}_2
\le\norm{Xz}_2
\le\chi_{\max}\norm{\wh L\wh Uz}_2.
\end{equation}
Let $y\in\mathbb{R}^n$ be arbitrary and let
$z=\wh R^{-1}y$.
Then~\eqref{eq:mplucqr-chi-vector-bound} becomes
\begin{equation*}
\chi_{\min}
\norm{\wh L\wh U\wh R^{-1}y}_2
\le
\norm{X\wh R^{-1}y}_2
\le
\chi_{\max}
\norm{\wh L\wh U\wh R^{-1}y}_2.
\end{equation*}
Applying \cref{lem:mplucqr-LUR-condition} to the two outer
terms gives
\begin{equation*}
0.965\cdot\chi_{\min}\norm{y}_2
\le
\norm{Yy}_2
\le
1.034\cdot\chi_{\max}\norm{y}_2.
\end{equation*}
Taking the minimum and maximum over all unit vectors $y$ proves
\eqref{eq:mplucqr-XR-bounds}.
\end{proof}

\begin{proof}[Proof of \cref{thm:mplucqr-convergence}]
As mentioned, it can be shown that the Cholesky factorization will not break down if \eqref{eq:mplucqr-assump-u-L} and~\eqref{eq:mplucqr-assump-ul} is satisfied, and the bound \eqref{eq:mplucqr-Rinvnorm} shows that $\wh R$ is nonsingular.

It remains to estimate the condition number of the computed preconditioned matrix $\wh Q$.
Multiplying $\wh Q\wh R-X=E_5$ on the right by $\wh R^{-1}$ gives
\begin{equation*}
\wh Q-Y=E_5\wh R^{-1}.
\end{equation*}
Using~\eqref{eq:mplucqr-E5-norm}, we obtain
\begin{equation}
\label{eq:mplucqr-QY-error}
\begin{aligned}
\norm{\wh Q-Y}_2
&\le\norm{E_5}_2\norm{\wh R^{-1}}_2\\
&\le n^2u\norm{\wh Q}_2 \norm{\wh R}_2\norm{\wh R^{-1}}_2\\
&=n^2u\cond_2(\wh R)\norm{\wh Q}_2.
\end{aligned}
\end{equation}
By~\eqref{eq:mplucqr-Rcond} and
assumption~\eqref{eq:mplucqr-assump-u-LU},
\begin{equation}
\label{eq:mplucqr-nuRcond}
n^2u\cond_2(\wh R)
\le
1.066n^2u\cond_2(\wh L)\cond_2(\wh U)
\le\frac{1.066}{64}
<0.0167.
\end{equation}

We first estimate the largest singular value of $\wh Q$.
By~\eqref{eq:mplucqr-QY-error}, \eqref{eq:mplucqr-nuRcond}, and~\eqref{eq:mplucqr-XR-bounds},
\begin{equation*}
\norm{\wh Q}_2\le\norm{Y}_2+\norm{\wh Q-Y}_2
\le1.034\cdot\chi_{\max}+0.0167\norm{\wh Q}_2.
\end{equation*}
Rearranging the inequality gives
\begin{equation}
\label{eq:mplucqr-Q-upper}
\norm{\wh Q}_2
\le
\frac{1.034}{1-0.0167}\cdot\chi_{\max}
<1.052\cdot\chi_{\max}.
\end{equation}

For the smallest singular value, by Weyl's theorem,
together with~\eqref{eq:mplucqr-XR-bounds} and~\eqref{eq:mplucqr-QY-error},
we have
\begin{equation*}
\begin{aligned}
\sigma_{\min}(\wh Q)
&\ge\sigma_{\min}(Y)-\norm{\wh Q-Y}_2\\
&\ge0.965\cdot\chi_{\min}
-n^2u\cond_2(\wh R)\norm{\wh Q}_2.
\end{aligned}
\end{equation*}
Using~\eqref{eq:mplucqr-Rcond} and~\eqref{eq:mplucqr-Q-upper} yields
\begin{equation*}
\begin{aligned}
\sigma_{\min}(\wh Q)
&>
0.965\chi_{\min}
-1.066\cdot1.052\,
 n^2u\cond_2(\wh L)\cond_2(\wh U)\chi_{\max}\\
&=
\left(
0.965
-1.066\cdot1.052\,
 n^2u\cond_2(\wh L)\cond_2(\wh U)
 \frac{\chi_{\max}}{\chi_{\min}}
\right)\chi_{\min}.
\end{aligned}
\end{equation*}
The second bound in assumption~\eqref{eq:mplucqr-assump-chi}
therefore gives
\begin{equation}
\label{eq:mplucqr-Q-lower}
\sigma_{\min}(\wh Q)
\ge\bigg(0.965-\frac{1.066\cdot1.052}{64}\bigg)\chi_{\min}
>0.947\chi_{\min}.
\end{equation}
Finally, combining~\eqref{eq:mplucqr-Q-upper} and
\eqref{eq:mplucqr-Q-lower}, we obtain
\begin{equation*}
\cond_2(\wh Q)=\frac{\norm{\wh Q}_2}{\sigma_{\min}(\wh Q)}<\frac{1.052}{0.947}\cdot\frac{\chi_{\max}}{\chi_{\min}}
\le1.12\bigl(1+(n\rho+\sqrt{n})\ul\cond_2(X)\bigr),
\end{equation*}
where the last inequality uses assumption~\eqref{eq:mplucqr-assump-chi}.
This proves~\eqref{eq:mplucqr-main-bound}.
\end{proof}

\bibliographystyle{plain-doi}
\bibliography{references}

@String{j-ACM                   = "J. ACM"}

@String{j-AN                    = "Acta Numerica"}

@String{j-ETNA                  = "Electron. Trans. Numer. Anal."}

@STRING{j-JRNBS                 = "J. Res. Natl. Bur. Stand."}

@String{j-MATH-COMP             = "Math. Comp."}

@string{j-NUM-ALG               = "Numer. Algorithms"}

@string{j-NUM-MATH              = "Numer. Math."}

@string{j-PAR-COMP              = "Parallel Comput."}

@String{j-SIMAX                 = "SIAM J. Matrix Anal. Appl."}

@String{j-SIREV                 = "SIAM Rev."}

@String{j-SISC                  = "SIAM J. Sci. Comput."}

@String{j-SISSC                 = "SIAM J. Sci. Statist. Comput."}

@String{pub-JH                  = "Johns Hopkins University Press"}

@String{pub-JH:adr              = "Baltimore, MD, USA"}

@String{pub-SIAM                = "Society for Industrial and Applied
                                   Mathematics"}

@String{pub-SIAM:adr            = "Philadelphia, PA, USA"}

@TechReport{bala22,
  author = "Balabanov, Oleg",
  title = "Randomized {Cholesky} {QR} factorizations",
  year = 2022,
  type = "{arXiv}:2210.09953 {[math.NA]}",
  month = feb,
  doi = "10.48550/arXiv.2210.09953"
  }

@Article{cahi17,
  author = "Erin Carson and Nicholas J. Higham",
  title = "A New Analysis of Iterative Refinement and its Application to
           Accurate Solution of Ill-Conditioned Sparse Linear Systems",
  journal = j-SISC,
  year = 2017,
  volume = 39,
  number = 6,
  pages = "A2834-A2856",
  doi = "10.1137/17M1122918"
  }

@Article{cahi18,
  author = "Erin Carson and Nicholas J. Higham",
  title = "Accelerating the Solution of Linear Systems by Iterative
           Refinement in Three Precisions",
  journal = j-SISC,
  year = 2018,
  volume = 40,
  number = 2,
  pages = "A817-A847",
  doi = "10.1137/17M1140819"
  }

@Article{clwo17,
  author = "Clarkson, Kenneth L. and Woodruff, David P.",
  title = "Low-Rank Approximation and Regression in Input Sparsity Time",
  journal = j-ACM,
  year = 2017,
  volume = 63,
  number = 6,
  pages = "54:1-54:45",
  doi = "10.1145/3019134"
  }

@Article{dsyg18,
  author = "Duersch, Jed A. and Shao, Meiyue and Yang, Chao and Gu, Ming",
  title = "A Robust and Efficient Implementation of {LOBPCG}",
  journal = j-SISC,
  year = 2018,
  volume = 40,
  number = 5,
  pages = "C655-C676",
  doi = "10.1137/17M1129830"
  }

@Article{fkny20,
  author = "Fukaya, Takeshi and Kannan, Ramaseshan and Nakatsukasa, Yuji and
                  Yamamoto, Yusaku and Yanagisawa, Yuka",
  title = "Shifted {Cholesky} {QR} for Computing the {QR} Factorization of
                  Ill-Conditioned Matrices",
  journal = j-SISC,
  year = 2020,
  volume = 42,
  number = 1,
  pages = "A477–A503",
  month = jan,
  doi = "10.1137/18m1218212"
  }

@InProceedings{fnyy14,
  author = "Fukaya, Takeshi and Nakatsukasa, Yuji and Yanagisawa, Yuka and
            Yamamoto, Yusaku",
  title = "{CholeskyQR2}: {A} Simple and Communication-Avoiding Algorithm For
           Computing a Tall-Skinny {QR} Factorization on a Large-Scale
           Parallel System",
  booktitle = "Proceedings of the 5th Workshop on Latest Advances in Scalable
               Algorithms for Large-Scale Systems",
  year = 2014,
  pages = "31-38",
  publisher = "IEEE",
  doi = "10.1109/ScalA.2014.11"
  }

@TechReport{gaip24,
  author = "Garrison, James E. and Ipsen, Ilse C. F.",
  title = "A Randomized Preconditioned {Cholesky}--{QR} Algorithm",
  year = 2024,
  month = jun,
  type = "{arXiv}:2406.11751 [math.NA]",
  doi = "10.48550/arXiv.2406.11751"
  }

@TechReport{gci26,
  author = "Garrison, James E. and Chen, Chao and Ipsen, Ilse C. F.",
  title = "Numerically stable {Cholesky}--{QR} on {GPU} via mixed-precision
           randomized preconditioning",
  year = 2026,
  type = "{arXiv}:2606.18411 {[math.NA]}",
  month = jun,
  doi = "10.48550/arXiv.2606.18411"
  }

@Article{goto08,
  author = "Gotsman, Craig and Toledo, Sivan",
  title = "On the Computation of Null Spaces of Sparse Rectangular Matrices",
  journal = j-SIMAX,
  year = 2008,
  volume = 30,
  number = 2,
  pages = "445-463",
  doi = "10.1137/050638369"
  }

@Book{gova13,
  author = "Golub, Gene H. and Van Loan, Charles F.",
  title = "Matrix Computations",
  publisher = pub-JH,
  year = 2013,
  address = pub-JH:adr,
  edition = "4th",
  isbn = "1-4214-0794-9"
  }

@TechReport{gzwc26,
  author = "Guan, Haoran and Zou, Zhenyu and Wei, Yufeng and Chen, Yipei and
            You, Peiting and Fan, Yuwei",
  title = "{RCLUPPr}: {A} New Randomized {CholeskyQR} with {LU} Preconditioning",
  year = 2026,
  type = "{arXiv}:2607.15561 {[math.NA]}",
  month = jul,
  doi = "10.48550/arXiv.2607.15561"
  }

@Book{high02,
  author = "Nicholas J. Higham",
  title = "Accuracy and Stability of Numerical Algorithms",
  publisher = pub-SIAM,
  year = 2002,
  address = pub-SIAM:adr,
  edition = "2nd",
  pages = "xxx+680",
  doi = "10.1137/1.9780898718027",
  isbn = "0-89871-521-0"
  }

@article{hima22,
  author = "Nicholas J. Higham and Theo Mary",
  title = "Mixed Precision Algorithms in Numerical Linear Algebra",
  journal = j-AN,
  volume = 31,
  pages = "347-414",
  month = may,
  year = 2022,
  doi = "10.1017/s0962492922000022",
  }

@Article{hmt11,
  author = "Halko, N. and Martinsson, P. G. and Tropp, J. A.",
  title = "Finding Structure with Randomness: {Probabilistic} Algorithms for
                  Constructing Approximate Matrix Decompositions",
  journal = j-SIREV,
  year = 2011,
  volume = 53,
  number = 2,
  pages = "217-288",
  month = jan,
  doi = "10.1137/090771806"
  }

@Article{hsby25,
  author = "Higgins, Andrew J. and Szyld, Daniel B. and Boman, Erik G. and
            Yamazaki, Ichitaro",
  title = "Analysis of Randomized {Householder--Cholesky} {QR} Factorization
           with Multisketching",
  journal = j-NUM-MATH,
  year = 2025,
  volume = 157,
  pages = "1695-1737",
  doi = "10.1007/s00211-025-01492-5"
  }

@Article{jeru13,
  author = "Jeannerod, Claude-Pierre and Rump, Siegfried M.",
  title = "Improved Error Bounds for Inner Products in Floating-Point Arithmetic",
  journal = j-SIMAX,
  year = 2013,
  volume = 34,
  number = 2,
  pages = "338–344",
  month = jan,
  doi = "10.1137/120894488"
  }

@Article{kms23,
  author = "Kressner, Daniel and Ma, Yuxin and Shao, Meiyue",
  title = "A mixed precision {LOBPCG} algorithm",
  journal = j-NUM-ALG,
  year = 2023,
  volume = 94,
  pages = "1653-1671",
  doi = "10.1007/s11075-023-01550-9"
  }

@Article{lanc50,
  author = "Lanczos, Cornelius",
  title = "An Iteration Method for the Solution of the Eigenvalue Problem Of
           Linear Differential and Integral Operators",
  journal = j-JRNBS,
  year = 1950,
  volume = 45,
  number = 4,
  pages = "255-282",
  doi = "10.6028/jres.045.026"
  }

@Article{matr20,
  author = "Martinsson, Per-Gunnar and Tropp, Joel A.",
  title = "Randomized Numerical Linear Algebra: {Foundations} and Algorithms",
  journal = j-AN,
  year = 2020,
  volume = 29,
  pages = "403-572",
  month = may,
  doi = "10.1017/s0962492920000021"
  }

@Article{neym02,
  author = "Neymeyr, Klaus",
  title = "A Geometric Theory for Preconditioned Inverse Iteration Applied To
           a Subspace",
  journal = j-MATH-COMP,
  year = 2002,
  volume = 71,
  number = 237,
  pages = "197-216",
  doi = "10.1090/S0025-5718-01-01357-6"
  }

@Article{ruje14,
  author = "Rump, Siegfried M. and Jeannerod, Claude-Pierre",
  title = "Improved Backward Error Bounds for {LU} and {Cholesky} Factorizations",
  journal = j-SIMAX,
  year = 2014,
  volume = 35,
  number = 2,
  pages = "684-698",
  doi = "10.1137/130927231"
  }

@Article{sasc86,
  author = "Saad, Youcef and Schultz, Martin H.",
  title = "{GMRES}: {A} Generalized Minimal Residual Algorithm for Solving
                  Nonsymmetric Linear Systems",
  journal = j-SISSC,
  year = 1986,
  volume = 7,
  number = 3,
  pages = "856–869",
  month = jul,
  doi = "10.1137/0907058"
  }

@Article{slvo96,
  author = "Sleijpen, Gerard L. G. and Van der Vorst, Henk A.",
  title = "A {Jacobi}--{Davidson} Iteration Method for Linear Eigenvalue Problems",
  journal = j-SIMAX,
  year = 1996,
  volume = 17,
  number = 2,
  pages = "401-425",
  doi = "10.1137/S0895479894270427"
  }

@Misc{tipr21,
  author = "Tisseur, Françoise and Pranesh, Srikara",
  title = "Mixed Precision {Cholesky}-{QR} Algorithm with Applications",
  howpublished = "Conference abstract, SIAM Conference on Computational Science
                  and Engineering (CSE21)",
  year = 2021,
  note = "March 3, 2021",
  url = "https://meetings.siam.org/sess/dsp_talk.cfm?p=107566"
  }

@Misc{tiss21,
  author = "Tisseur, Françoise",
  title = "Mixed Precision Tall and Thin {QR} Factorization with Applications",
  howpublished = "Conference abstract, New Directions in Numerical Linear
                  Algebra and High Performance Computing: Celebrating the 70th
                  Birthday of Jack Dongarra",
  year = 2021,
  note = "July 7--8, 2021, p.~3. Joint work with Srikara Pranesh",
  url = "https://nlagrouporg.wordpress.com/wp-content/uploads/2021/06/talk-abstracts.pdf"
  }

@Article{too20,
  author = "Terao, Takeshi and Ozaki, Katsuhisa and Ogita, Takeshi",
  title = "{LU}-{Cholesky} {QR} Algorithms for Thin {QR} Decomposition",
  journal = j-PAR-COMP,
  year = 2020,
  volume = 92,
  pages = 102571,
  month = apr,
  doi = "10.1016/j.parco.2019.102571"
  }

@Article{ynyf15,
  author = "Yusaku Yamamoto and Yuji Nakatsukasa and Yuka Yanagisawa and Takeshi
                  Fukaya",
  title = "Roundoff Error Analysis of the {CholeskyQR2} Algorithm",
  journal = j-ETNA,
  year = 2015,
  volume = 44,
  pages = "306-326",
  url = "https://etna.ricam.oeaw.ac.at/volumes/2011-2020/vol44/abstract.php?vol=44&pages=306-326"
  }

@Article{ytd15,
  author = "Yamazaki, Ichitaro and Tomov, Stanimire and Dongarra, Jack",
  title = "Mixed-Precision {Cholesky} {QR} Factorization and Its Case Studies on
                  Multicore {CPU} with Multiple {GPUs}",
  journal = j-SISC,
  year = 2015,
  volume = 37,
  number = 3,
  pages = "C307–C330",
  month = jan,
  doi = "10.1137/14m0973773"
  }

\end{document}